\documentclass[12pt]{article}
\usepackage[dvips]{graphics}
\usepackage{amsmath}
\usepackage{amssymb}
\usepackage{cite}
\usepackage{color}
\usepackage{epsfig} 
\def\X{\mathbf{x}}
\def\Z{\mathbf{z}}
\def\W{\mathbf{w}}
\usepackage{amsthm}
\usepackage{hyperref}

\newtheorem{prop}{Proposition}

\newtheorem{example}{Example}

\title{\Large Homeostasis phenomenon in predictive inference when using a wrong learning model: a tale of random
	split of data into training and test sets\footnote{The research note is prepared while writing a discussion for Professor Efron's paper \cite{Efron2020}. It contains a more detailed and elaborated discussion, including additional technical details. 
The research is supported in part by NSF-DMS1737857 and DMS1812048.
}}
	
\author{Min-ge Xie and Zheshi Zheng \\ Department of Statistics, Rutgers University}
\date{February 2020}

\begin{document}

\maketitle

This note uses a conformal prediction procedure to provide further support on several points discussed by Professor Efron \cite{Efron2020} concerning 
prediction, estimation and IID assumption. 
It aims to convey the following messages:  
\begin{itemize} 
\item  Under the IID (e.g., random split of training and testing data sets) assumption, prediction is indeed an easier task than estimation, since prediction has a {\it homeostasis} property in this case --- Even if the model used for learning is completely wrong, 
the prediction results maintain valid.

\item If the IID assumption is violated  (e.g., a targeted prediction on specific individuals), the homeostasis property is often disrupted  and the prediction results under a wrong model are usually invalid.
\item
Better model estimation typically leads to more accurate prediction~in~both IID and non-IID cases. Good modeling and  estimation practices are important and, in many times, crucial for obtaining  good prediction~results. 
\end{itemize}
The discussion also provides one explanation why the deep learning method works so well in academic exercises (with experiments set up by randomly splitting the entire data into training and testing data sets), but fails to deliver many ``killer applications'' in real world applications.
\section{Introduction}

This outstanding paper by Professor Efron \cite{Efron2020} provides stimulating discussions on the future of our field in the remainder of the 21st Century. 
In this note, we echo and also provide additional support to two important points made by Professor Efron: (1) prediction is ``an easier task than either
attribution or estimation''; (2) the IID assumption (on both training and testing  data sets) is
crucial
in the current developments on predictions, but we also need to do more for the case when the IID assumption is not met. 
Based on our own research, we  provide additional evidence to support these discussions, including a discovery why prediction has a {\it homeostasis} property and works well under the IID setting even if the learning model used is completely wrong. We specifically highlight the importance of having a
good learning model with good  estimation to obtain a good prediction.  We provide examples to show that, for the task of prediction, a good modeling  and inference practice is important in the IID case and it becomes essential
for non-IID case. 
The message remains: to get a good prediction outcome, we still need to make effort to build a good learning model and estimation algorithm, even if sometimes prediction appears to be an easier task than estimation. 

From the outset, we would like to comment on that it is not a straw-man argument to consider non-IID
testing data. On the contrary, they are prevalent in data science. In addition to those examples provided by Professor Efron~that showed ``drift,'' we can easily imagine non-IID examples in many typical applications. For instance, a predictive algorithm is trained on 
a 
database of patient medical records and we would like to predict  potential outcomes of a treatment for a new patient with  more severe symptoms than what the average patient shows. The new patient with more severe symptoms is not a typical IID  draw from the general patient population. Similarly, in the finance sector, one is often interested in predicting the financial performance of a particular company or group. If a predictive model is trained on all institutes, then the testing data~(of the specific group of companies of interest) are unlikely IID draws from the same general population of the learning data. The limitation of IID assumption has hampered our efforts to fully take advantage of fast-developing machine learning methodologies (e.g., deep neural network model, tree based methods, etc.) in many real-world applications,  a point that we will have more elaboration later. 

Our discussions in this note are based on a so-called {\it conformal prediction} procedure, an attractive new prediction framework that is error (or model) distribution
free;
cf., e.g., \cite{vovk2005algorithmic,shafer2008tutorial,lei2018distribution,barber2019predictive, barber2019limits}.
We demonstrate that, under the IID assumption, the predictive conclusion is always valid even if the model used to train the data is completely wrong. We discover a homeostasis phenomenon that the prediction is resistant to wrong learning models in the IID case because 
the expected bias caused by learning using a wrong model can largely be offset by the corresponding negatively shifted predictive errors (cf., Sections~\ref{sec:2.3} and~\ref{sec:3.1}).
This robustness result clearly supports the claim that prediction is an easier task than modeling and estimation. However, the use of a wrong learning model has at least two undesirable impacts on prediction: (a) A prediction based on a wrong model typically produces a much wider predictive interval (or a wider predictive curve) than that based on a correct model; 
(b) 
Although the IID case enjoys a nice homeostatic cancellation of bias (in fitted model) and shifts (in associated predictive errors) when using a wrong learning model,
in the non-IID case this cancellation is often no longer effective, resulting in invalid predictions. 
The use of a correct learning 
model can help mitigate and sometimes solve the problem of invalid  prediction for non-IID (e.g., drifted or individual-specific) testing~data.

The rest of the note is arranged as follows. Section 2 describes the conformal predictive inference in general terms. The prediction is valid under the IID setting, even if the learning model used is completely wrong.
Section~3 contains two case studies, one on linear regression and the other on neural network model, to study the impact of using a wrong learning model 
on prediction under both IID and non-IID settings.  
Concluding remarks~are~in~Section~4.

\section{Prediction, testing data and learning models}

As in equation (6.4) of Professor Efron's article, we consider a typical setup in data science: Suppose we have a training (observed) data set of size $n$:
${\cal D}_{obs} = \{({\X}_i, y_i), i = 1,\ldots, n\}$, where $({\X}_i, y_i)$, $i = 1,\ldots, n$, are IID random samples from an unknown population ${\cal F}$. For a given ${\X}_{new}$, we would like to predict what $y_{new}$ would be. 
We first use the typical assumption that $({\X}_{new}, y_{new})$ is also an IID draw from ${\cal F}$. Later we relax this requirement and only assume that $y_{new}|{\X}_{new}$ relates to $\X_{new}$ the same way as $y_i|{\X}_i$ relates to $\X_i$, 
 but ${\X}_{new}$ is fixed or follows a marginal distribution that is different from that of ${\X}_i$. 

For notation convenience, we consider $({\X}_{new}, y_{new})$ as the $(n+1)$-th observation and introduce the index $n+1$, with $\X_{n+1} = \X_{new}$ and $y_{n+1}$ as a potential value (or a ``guess'') of the unobserved $y_{new}$.  Unless specified otherwise, the index ``$n + 1$'' and index ``new'' are exchangeable throughout the note.  

\subsection{Conformal prediction and level $(1 - \alpha)$ conformal predictive intervals} A {\it conformal prediction} procedure is a distribution free prediction method that
has attracted increasing attention in computer science and statistical learning communities in recent years; cf., e.g.,  \cite{vovk2005algorithmic,shafer2008tutorial,lei2018distribution,barber2019predictive, barber2019limits}. 
The idea of conformal prediction is straightforward. In order to make a prediction of the unknown $y_{new}$ given $\X_{n+1} = \X_{new}$, we examine a potential value $y_{n+1}$, and see how ``conformal" the pair $({\X}_{n+1}, y_{n+1})$ is among the observed $n$ pairs of IID data points $({\X}_i, y_i)$, $i = 1, \ldots, n$. The higher the ``conformality,'' the more likely $y_{new}$ takes the potential value $y_{n+1}$. Frequently, a learning model, say $y_i \sim \mu({\X}_i)$ for $i = 1 \ldots, n, n+1$, is used to assist prediction. However, the learning model is not essential. 
As we will see later, even if $\mu(\cdot)$ is totally wrong or does not exist, a conformal prediction can still provide us valid prediction, as long as the IID assumption holds for both the training and testing data, i.e., $(\X_i, y_i), (\X_{new}, y_{new}) \overset{\tiny iid}{\sim} {\cal F}$, for $i = 1, \ldots, n$.

To be specific, this note employs a conformal prediction procedure that is referred to as the {\it Jackknife-plus method} by \cite{barber2019predictive}.
Specifically, consider a combined collection of both the training and testing data but with the unknown $y_{new}$ replaced by a potential value $y_{n+1}$: ${\cal A} = {\cal D}_{obs} \cup \{(\X_{new}, y_{n+1})\} = \{({\X}_{i}, y_i), i = 1, \ldots, n, n+1\}$. 
We define {\it conformal residuals} 
$R_{ij}(y) = y_{i} -  \hat y_{i}^{- (i, j)},$ for $i \not = j \,\hbox{and}\, i,j = 1,\ldots, n, n+1$, where $\hat y_{i}^{- (i, j)}$ 
is the prediction of $y_{i}$ based on the leave-two-out dataset ${\cal A}^{-(i,j)} = {\cal A} - \{({\X}_i, y_i), ({\X}_j, y_j)\}$.
If a working model $\mu(\cdot)$ is used, for instance, the model is first fit based on the leave-two-out dataset ${\cal A}^{-(i,j)}$ and the point prediction is set to be $\hat y_{i}^{- (i, j)} = \hat \mu(\X_i; {\cal A}^{-(i,j)})$, 
where $\hat \mu(\cdot; {\cal A}^{-(i,j)})$ is the fitted (trained) model using ${\cal A}^{-(i,j)}$. 

For each given $y_{n+1}$ (a potential value of $y_{new}$),  we define \begin{equation}
\label{eq:Q}
    Q_{n}(y_{n+1}) = \frac{1}{n} \sum_{i = 1}^n {\bf 1}_{\{R_{n+1,i} \ge R_{i, n+1}\}},
\end{equation} 
which relates to the degree of 
``conformity'' of the residual values $R_{n+1,i} = y_{n+1} - \hat y_{n+1}^{-(i, n+1)}$ among the residuals (which in fact are the leave-one-out residuals of using only the training data ${\cal D}_{obs}$) $R_{i, n+1} =  y_{i} -  \hat y_{i}^{- (i, n+1)}$, $i = 1, \ldots, n$. If $Q_{n}(y_{n+1}) \approx \frac 12$, then $R_{n+1,i}$ is around the middle of the training data residuals $R_{i, n+1}$ and thus ``most conformal.'' When $Q_{n}(y_{n+1}) \approx 0$ or $\approx 1$,  $R_{n+1,i}$ is at the extreme ends of the training data residuals $R_{i, n+1}$ and thus ``least conformal.'' This intuition leads us to define a conformal predictive interval of $y_{new}$ as
\begin{align}\label{eq:PI}
C_{\alpha} & =
\left\{y: Q_{n}(y) \ge \frac{\alpha}{2}\right\} \bigcap \left\{y: 1- Q_{n}(y) \ge \frac{\alpha}{2}\right\} \nonumber \\
 & = 
\left[q_{\frac{\alpha}{2}}\left(\{\hat y_{n+1}^{-(i, {n+1})} +   R_{i,n+1}\}_{i =1}^n\right), \, q_{1 - \frac{\alpha}{2}}\left(\{\hat y_{n+1}^{-(i, {n+1})} + R_{i,n+1}\}_{i =1}^n\right)\right],
\end{align}
where $q_{\alpha}\left(\{a_i\}_{i =1}^n\right)$ is the $\alpha$-th quantile of $a_1, \ldots, a_n$. 
The predictive interval (\ref{eq:PI}) is a slightly variant version of the {\it Jackknife-plus predictive interval} proposed by \cite{barber2019predictive}.

The following proposition states that, under the IID assumption,  $C_\alpha$ defined in~(\ref{eq:PI})
is guaranteed a level-$(1 - 2\alpha)$ predictive set for $y_{new}$. 
We outline a proof of the 
proposition in Supplementary.  The proposition and proof is almost the same as that provided in Theorem~1 of  \cite{barber2019predictive}, except that the absolute residuals $|R_{ij}| = |y_{i} -  \hat y_{i}^{- (i, j)}|$ are used instead throughout their development. 

\begin{prop}
\label{LM-1}
Suppose $(\X_i, y_i), (\X_{new}, y_{new}) \overset{\tiny iid}{\sim} {\cal F}$, for $i = 1, \ldots, n$. Then, we have $\mathbb{P}\left(y_{new} \in C_\alpha \right) \ge 1 - 2 \alpha$. 
\end{prop}

\noindent
Proposition~\ref{LM-1} is proved for a finite $n$, with a (conservatively) guaranteed coverage rate of  $(1-2 \alpha)$.~\cite{barber2019predictive} pointed out empirically $C_\alpha$ has a typical coverage rate of $1 - \alpha$.  In the rest of the note, we treat $C_{\alpha}$ as an approximate level-$(1 - \alpha)$ predictive interval. 

A striking result is that Proposition~\ref{LM-1} holds, even if the learning model $\mu(\cdot)$ used to assist prediction is completely wrong, as long as $\hat \mu(\,\cdot\,; {\cal A}^{-(i,j)})$ and $\hat \mu(\,\cdot\,; {\cal A}^{-(i',j')})$, for any two pairs $(i,j)$ and $(i',j')$, $i \not =j$, $i'\not = j'$, maintain  ``symmetry'' or
``exchangeability'' (when shuffling indices)  due to the IID assumption. This amazingly robust property is highly touted in the machine learning community. It gives support to the sentiment of using  ``black box'' algorithms where the role of model fitting is reduced to an afterthought, although we will provide arguments to counter this sentiment later in the note. 

\subsection{Conformal predictive distribution and  predictive curve} To get a full picture of the prediction intervals at all significance levels (as we present later in Figures \ref{fig:PV_LM} and \ref{fig:PV-NN}), 
we would like to briefly describe the notions of {\it predictive distribution} (cf.,  \cite{Lawless2005,Shen2018,Vovk2019}) and {\it predictive curve}. Predictive distribution in Bayesian inference is well known, but the development of predictive distribution with confidence interpretation is relatively new; cf., \cite{Lawless2005,Shen2018}.
Note that a 
predictive interval 
has the same frequency interpretation as a
confidence interval, except that it is developed for a random $y_{new}$ instead of a parameter of interest. Similarly, a  predictive distribution (with a  confidence interpretation) can be viewed as an extension of a confidence distribution but developed for the random $y_{new}$ instead for a parameter of interest. 

\cite{cox1958} suggested that a  confidence distribution be introduced ``in terms of the set of confidence intervals of all levels".
To better understand the concept of predictive distribution and predictive curve, especially how to relate them to predictive intervals of all levels, it is prudent to briefly take a look at confidence distribution and confidence curve, and then move on to prediction. 
We consider 
a toy example below.  

\begin{example}
Assume in this toy example that  $y_1, \ldots, y_n \overset{\tiny iid}{\sim} N(\theta, 1)$. Instead of using a point ($\bar y = \frac1n \sum_{i = 1}^n y_i$) or an interval ( $\bar y \pm \frac{1}{\sqrt{n}}\Phi^{-1}(1 - \frac{\alpha}2)$)
to  estimate the unknown  parameter $\theta$, 
a confidence distribution  suggests to use a sample-dependent function $N(\bar y, \frac1n)$, or more formally in the cumulative distribution function form $H_n(\theta) = \Phi(\sqrt{n}(\theta - \bar y))$,  to estimate the unknown parameter $\theta$; cf, e.g., \cite{Efron1993, Xie2013, Schweder2016}. A nice feature of a confidence distribution is that it can represent 
confidence intervals of all levels. 
For example, the level-$(1 - \alpha)$ one-sided interval 
$(-\infty, \bar y + \frac{1}{\sqrt n} \Phi^{-1}(1 - \alpha)) = (-\infty, H_n^{-1}(1 -\alpha))$ and the level-$(1 - \alpha)$ two-sided interval $(\bar y + \frac{1}{\sqrt n} \Phi^{-1}(\frac{\alpha}2),\bar y + \frac{1}{\sqrt n}\Phi^{-1}(1 - \frac{\alpha}2)) = ( H_n^{-1}(\frac{\alpha}2), H_n^{-1}(1- \frac{\alpha}2))$. Here,  $H_n^{-1}(\cdot)$ is the inverse function of $H_n(\cdot)$. 

A closely related concept is  confidence curve $$
CV_n(\theta) = 2\min\{H_n(\theta), 1 - H_n(\theta)\},
$$ 
which was first introduced by \cite{Birnbaum1961} as an ``omnibus
form of estimation” that “incorporates confidence limits and intervals at all levels.''
For any $\alpha \in (0,1)$,  $\{\theta: CV_n(\theta) \ge \alpha\} $ is a level-$(1 -\alpha)$ two-sided confidence interval. We could view
the function $CV_n(\theta)$ as a result of stacking up  two-sided confidence intervals of all levels $1 - \alpha$ for $\alpha$ going from $0$ to $1$; cf., Figure 1 (a). The plot of confidence curve function $CV_n(\theta) = 2 \min\{\Phi(\sqrt{n}(\theta - \bar y)), 1 - \Phi(\sqrt{n}(\theta - \bar y)) \}$ provides a full picture of confidence intervals of all levels $1- \alpha \in (0,1)$,  with a peak point corresponding to a median unbiased estimator $\hat \theta_M = \bar y$ with $\mathbb{P}(\hat \theta_M \leq \theta) \geq \frac12$ and $\mathbb{P}(\hat \theta_M \geq \theta) \geq \frac12$. 
\begin{figure}
    \centering
    \includegraphics
    [width=5.6cm,height=3.8cm]{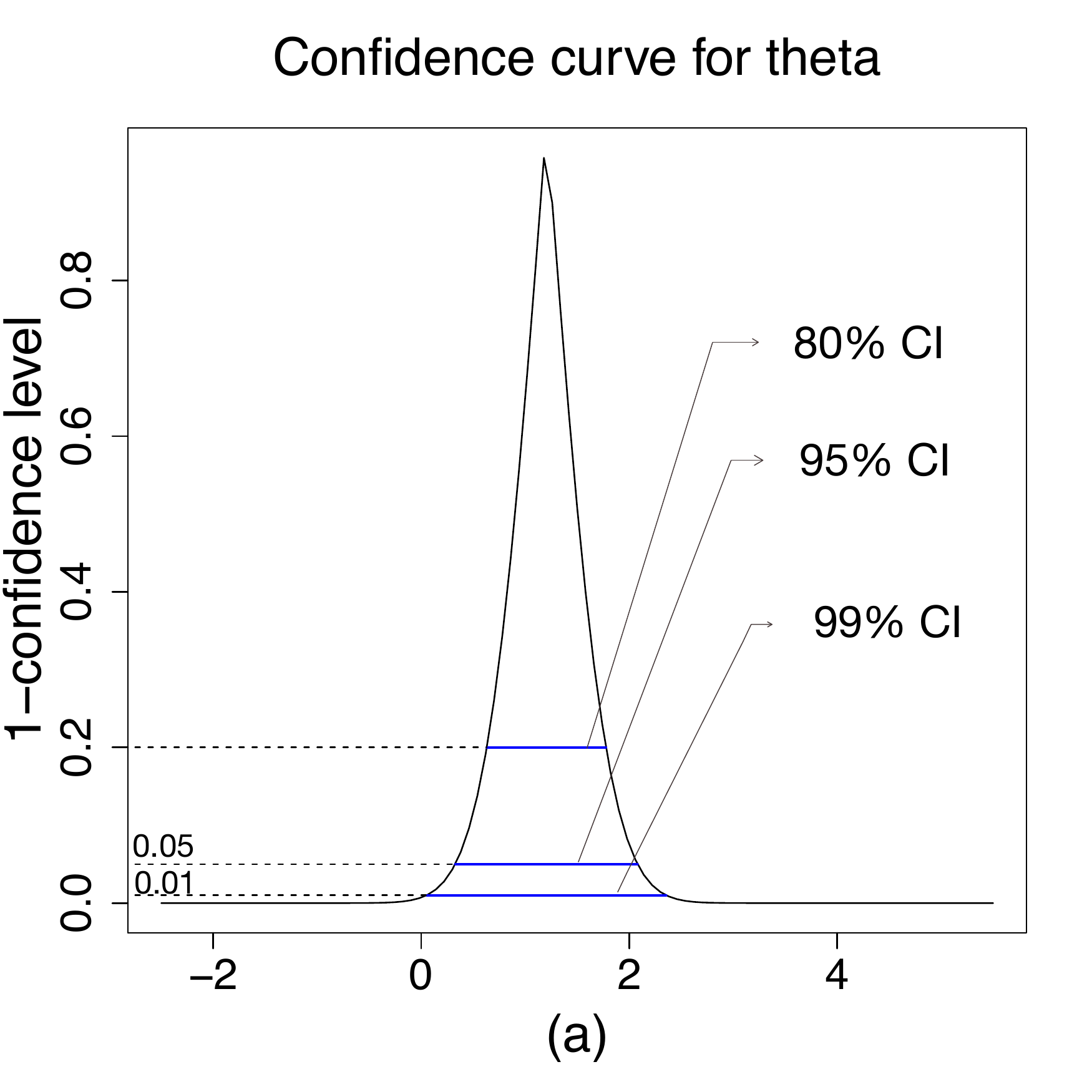}
    \includegraphics
    [width=5.6cm,height=3.8cm]{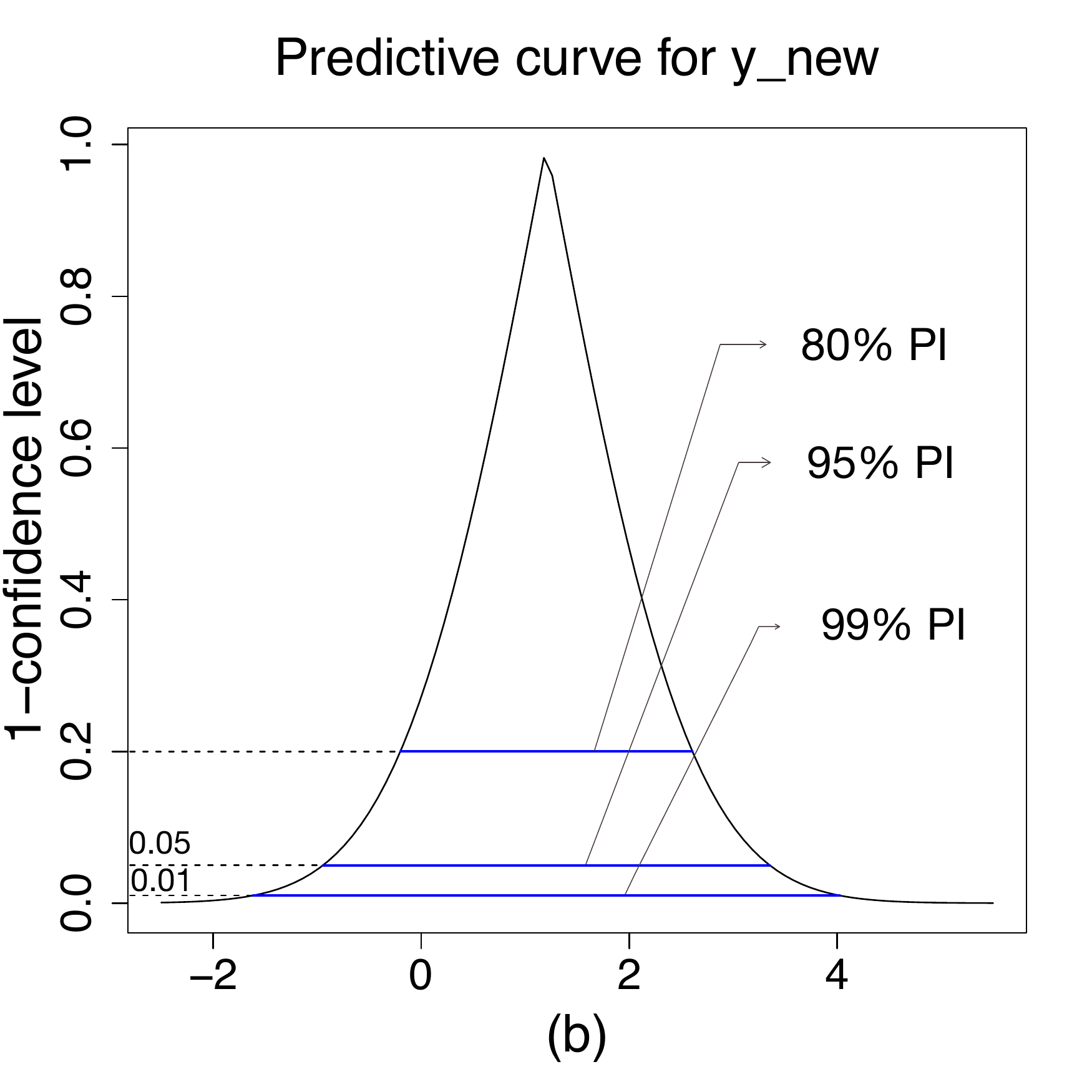}
    \caption{\small
    Plot of (a) confidence curve function $CV_n(\theta) = 2 \min\{\Phi(\frac{\theta - \bar y}{\sqrt{n}}), 1 - \Phi(\frac{\theta - \bar y}{\sqrt{n}}) \}$; (b) predictive curve function $PV_n(y) = 2\min\{\Phi\big(\frac{y - \bar y}{\sqrt{1+1/n}} \big), 1 - \Phi\big(\frac{y - \bar y}{\sqrt{1+1/n}} \big)\}$. 
    The plots provide a full picture of (a) confidence intervals and (b) predictive intervals of all levels. In particular, the curves can be formed as stacking up the endpoints of (a) the confidence intervals or (b) the predictive intervals at all levels of $1-\alpha$ for $\alpha$ from $0$ to $1$. The peak point corresponds to the median unbiased (a) point estimator $\hat \theta_M$ of $\theta$ and (b) point prediction $\hat y_M$ of $y_{new}$, respectively. The sample data used to generate the plots are from $N(1.35, 1)$ with $n = 5$.    
    }
    \label{fig:1}
\end{figure}

For a new sample $y_{new} \sim N(\theta, 1)$, a predictive distribution is $N(\bar y, 1 + \frac1n)$, or in its cumulative distribution function form
$Q_n(y) = \Phi\big(\frac{y - \bar y}{\sqrt{1+1/n}} \big)$. Parallel to confidence curve, we can define a predictive curve
{\small 
\begin{equation}
 \label{eq:PVa}
PV_n(y)  = 2\min\{Q_n(y), 1 - Q_n(y)\}  = 
2\min\{\Phi\big(\frac{y - \bar y}{\sqrt{1+1/n}} \big), 1 - \Phi\big(\frac{y - \bar y}{\sqrt{1+1/n}} \big)\}.
\end{equation}
}

\vspace{-3mm}
\noindent
Figure~\ref{fig:1} (b) is a plot of the predictive curve 
in (\ref{eq:PVa}). Again, we can view
the function $PV_n(y)$  as a result of stacking up two-sided predictive intervals of all levels $1 - \alpha$ for $\alpha$ going from $0$ to $1$. 
The plot of the predictive curve $PV_n(y) =$ {\small $ 2\min\big\{\Phi\big(\frac{y - \bar y}{\sqrt{1+1/n}} \big),$ $1 - \Phi\big(\frac{y - \bar y}{\sqrt{1+1/n}} \big)\big\}$} provides a full picture of predictive intervals~of all levels {\small $1- \alpha \in$} $(0,1)$. 
The peak point in Figure~\ref{fig:1}(b) corresponds to a median~un-
biased point predictor $\hat y_M = \bar y$ with $\mathbb{P}(\hat y_M \leq y_{new}) \geq \frac12$ and~$\mathbb{P}(\hat y_M \geq y_{new}) \geq \frac12$.
\end{example}

Back to our conformal prediction development, the function $Q_n(y)$ defined in (\ref{eq:Q}) is in essence a predictive distribution of $y_{new}$. 
The associated 
{\it predictive curve} for $y_{new}$ can then be defined as
$$PV_n(y) = 2 \min\{Q_{n}(y), 1 - Q_{n}(y)\}.$$ 
The 
predictive interval in (\ref{eq:PI}) is $C_\alpha = \{y: PV_n(y) \ge \alpha\}$.  We later plot our predictive curves $PV_n(y)$ in Figures~\ref{fig:PV_LM} and \ref{fig:PV-NN}, which provides a full picture of conformal predictive intervals of all levels in various setups.

Note that, in Example 1, $H(t)$ is the $p$-value for the one-sided test $H_0: \theta \leq t$ versus $H_1: \theta > t $ and $CV(t)$ is the is the $p$-value for the two-sided test $H_0: \theta = t$ versus $H_1: \theta \not = t$; cf., e.g., \cite{Xie2013,Schweder2016}. Thus, $H(\theta)$ and $CV(\theta)$ can be interpreted as the same quantities of $p$-value functions of one-sided and two-sided tests, respectively. Similarly, the predictive function $Q_n(y)$ and predictive curve $PV(y)$ also have the corresponding interpretation of $p$-value functions of one-sided test $H_0: y_{new} \leq y$ versus $H_1: y_{new} > y$ and two-sided test $H_0: y_{new} = y$ versus $H_1: y_{new} \not = y$, respectively.

\subsection{Validity vs efficiency and IID vs non-IID under a wrong learning model}\label{sec:2.3}

Although the validity result in Proposition~\ref{LM-1} is robust against wrong learning models under the IID setting, there is no free lunch.  
The predictive intervals obtained under a wrong model will typically be wider. For instance,
suppose that the true model is $y = \mu_0(\X) + \epsilon$, but a wrong model $y = \mu_1(\X) + e$ is used. Since $y = \mu_0(\X) + \epsilon =  \mu_1(\X) + \{\mu_0(\X) - \mu_1(\X)\} + \epsilon$, we have $e =  \{\mu_0(\X) - \mu_1(\X)\} + \epsilon$. 
So, when $\epsilon$ is independent of $\X$, ${\rm var}(e) =  {\rm var}(\{\mu_0(\X) - \mu_1(\X)\}) + {\rm var}(\epsilon) \ge {\rm var}(\epsilon)$ and the equality holds only when $\mu_1(\X) = \mu_0(\X)$. Thus, the error term $e$ under a wrong model has a larger variance than that of the error term $\epsilon$ under the true model. 
The larger ${\rm var}(\{\mu_0(\X) - \mu_1(\X)\})$ are (i.e. the more discrepant  $\mu_1(\X) - \mu_0(\X)$ are),
the larger the variance of the error term $e$ are. 
A larger error typically translates to less accurate estimation and prediction.

We have an intuitive explanation 
why a conformal predictive algorithm can still provide valid prediction even under a totally wrong model in the IID setting. Specifically, when we use a wrong model $\mu_1(\X)$, the corresponding point predictor will be biased by the magnitude of $\mu_1(\X_{new}) - \mu_0(\X_{new})$, but at the same time the error term $e$ absorbs the {\it bias}, thus producing residuals with a {\it shift} by the magnitude of $\mu_0(\X_i) - \mu_1(\X_i) = -\{\mu_1(\X_i) - \mu_0(\X_i)\}$. In a conformal prediction algorithm, the quantiles of residuals are added back to the point prediction to form the bounds of predictive intervals. If the IID assumption holds, the bias is offset by the shift. Along with greater  residual variance, the offsetting
helps ensure the validity of the conformal prediction.  
We call this tendency of self balance to maintain validity a {\it homeostasis phenomenon}, and will 
explain it in explicit mathematical terms
under a linear model in
Section~\ref{sec:3.1}.

The IID assumption is a 
crucial 
condition to ensure the validity of a prediction under a wrong model. If the IID assumption does not hold for the testing data, the prediction based on a wrong model (or a correct model but a wrong parameter estimation) is often invalid with huge errors, as we see in our case studies in Section~\ref{sec:3}. We think this IID assumption
also explains why deep neural networks and other machine learning methods work so well in academic research settings (where random split of data into training and testing sets is a common practice) but fail to produce ``killer applications'' to make predictions for a given patient or company whose $\X_{new}$ are often not close to the center of the training data. 
The good news is that, if we use a correct model for training and can get good model estimates, a reasonably acceptable prediction for a fixed  $\X_{new}$ is possible. This is illustrated in the case studies in Section~\ref{sec:3} under both linear and neural  network models. Indeed, modeling and estimation 
remain relevant and often crucial for prediction in both IID and non-IID cases.

\section{Case studies: conformal prediction under specific models}\label{sec:3}
\subsection{Prediction with data from a linear regression model}\label{sec:3.1} 
We assume in this subsection that the training data $({\X}_i, y_i)$, $i = 1, \ldots, n$, and $({\X}_{new}, y_{new})$ are from the following linear model:
\begin{equation}\label{eq:LM-mu0}
    y_i = \mu_0
    (\X_i) + \epsilon_i = {\X}_i^T\beta + \epsilon_i
\end{equation}
where $\beta$ is the unknown regression coefficient and $\epsilon_i$ are IID random errors with mean $0$ and variance $\sigma^2$. 
We would like to compare the performances of the Jacknife-plus prediction procedure 
under two different learning models:
$$
\hbox{(a) the true model: $\mu_0({\X}_i) = {\X}_i^T\beta$ \, vs \, (b) a wrong model: $\mu_1({\X}_i) = {\Z}_i^T\gamma$ 
},
$$
where ${\Z}_i$ is the first $q$ elements of the $p$ covariates of ${\X}_i$, $q < p$, and $\gamma$ is the corresponding $q \times 1$ unknown regression coefficient. We define notations: ${\bf Y}$ is the $n \times 1$ response vector of the training (observed) data, $\mathbf{X}$ and $\mathbf{Z}$ are the $n \times p$ and $n \times q$ design matrices, respectively, and we have a matrix partition $\mathbf{X} = (\mathbf{Z}, \,\, \mathbf{W})$.

Under the true learning model $\mu_0({\X}_i)$
and from the least squares estimation,  
we have, 
for each given $i$ and $s = i$ or $n+1 \, (new)$, 
\begin{align*}
    \hat y_s^{-(i,n+1)} & = \hat\mu_0(\X_s; {\cal A}^{-(i,n+1)} )  = \X_{s}^T(\mathbf{X}^T\mathbf{X} - \X_i \X_i^T)^{-1} (\mathbf{X}^T {\bf Y} - \X_i y_i) \\ \nonumber
    & = ... =
    \X_{s}^T \hat \beta - (y_i - \X_i^T \hat \beta) \frac{h_{is}}{1 - h_{ii}}, 
\end{align*}
where $\hat \beta = (\mathbf{X}^T\mathbf{X})^{-1} \mathbf{X}^T{\bf Y}$ is the least squares estimator using all training data (of size $n$) and $h_{ii} = \X_i^T  (\mathbf{X}^T \mathbf{X})^{-1}\X_i$ and $h_{is} =\X_s^T  (\mathbf{X}^T \mathbf{X})^{-1} \X_i$.
Therefore, for $i = 1, \ldots, n$ and replacing index $n+1$ with index $new$, 
\begin{align*}
 & \hat y_{n+1}^{-(i,{n+1})} +   R_{i,n+1} = \hat y_{new}^{-(i, new)} +   y_i - \hat y_{i}^{-(i,{new})} \\
& \qquad = \X_{new}^T \hat \beta + (y_i - \X_i^T \hat \beta)  \frac{1-h_{i, new}}{1 - h_{ii}} 
= \X_{new}^T \hat \beta + (1 - h_{i, new}) u_i,
\end{align*}
where $u_i = \frac{y_i-\X_i^T\hat\beta}{1 - h_{ii}}$ is the deleted residual (using all training data of size $n$) and $h_{i,new}$ $=\X_{new}^T  (\mathbf{X}^T \mathbf{X})^{-1} \X_i$. 
Thus, from (\ref{eq:PI}), the predictive interval of $y_{new}$ is: 
\begin{equation}
\label{eq:PI-true}
 [\X_{new}^T \hat \beta + q_{\frac{\alpha}{2}}(\{(1 - h_{i, new}) u_i\}_{i=1}^n), \,  \X_{new}^T \hat \beta + q_{1 - \frac{\alpha}{2}}(\{(1 - h_{i, new}) u_i\}_{i=1}^n)].  
\end{equation}
Note that, given $\X_{new}$, the point predictor $\X_{new}^T \hat \beta$ is an unbiased estimator of $E(y_{new}|$ $\X_{new}) = \X_{new}^T \beta$ and  $E\{(1 - h_{i, new}) u_i | \X_{new}\} = 0$, for $i = 1, \ldots, n$. 
Thus, the prediction interval (\ref{eq:PI-true}) can be interpreted as an interval ``centered'' at the unbiased predictor $\X_{new}^T\hat \beta$ with its width determined by the ``spread'' of the mean-zero ``noises'' $\{(1 - h_{i, new}) u_i, i = 1, \ldots, n\}$.   

When the wrong model $\mu_1(\Z)$ is used, we can use a similar derivation to get the predictive interval of $y_{new}$:
\begin{equation}
\label{eq:PI-wrong}
    \left[\Z_{new}^T \hat \gamma + q_{\frac{\alpha}{2}}(\{(1 - g_{i, new}) v_i\}_{i=1}^n), \, \Z_{new}^T \hat \gamma + q_{1 - \frac{\alpha}{2}}(\{(1 - g_{i, new}) v_i\}_{i=1}^n)\right].  
\end{equation}
where $\hat \gamma = (\mathbf{Z}^T \mathbf{Z})^{-1}\mathbf{Z}^T{\bf Y}$ is the least squares estimator using the wrong model, $g_{ii} = \Z_i^T  (\mathbf{Z}^T \mathbf{Z})^{-1}\Z_i$, $g_{i,new} =\Z_{new}^T  (\mathbf{Z}^T \mathbf{Z})^{-1} \Z_i$ and $v_{i} = \frac{y_i-\Z_i^T\hat \gamma}{1-g_{ii}}$. Here, given $\X_{new}$, the point predictor $\Z_{new}^T \hat \gamma$ is actually biased, with mean $\mathbb{E} \left(\Z_{new}^T\hat\gamma |\mathbf{X},\X_{new}\right)$  $=  \Z_{new}^T(\mathbf{Z}^T\mathbf{Z})^{-1}\mathbf{Z}^T \mathbf{X} \beta$. 
Thus, the bias caused by missing the covariates $\W_i$ is
\begin{align}
    bias &= \mathbb{E} \left[\Z_{new}^T\hat\gamma |\mathbf{X},\X_{new}\right] -\X_{new}^T\beta
    = - \mathbf{w}_{new}^T \beta_2 + \Z_{new}^T(\mathbf{Z}^T\mathbf{Z})^{-1}\mathbf{Z}^T\mathbf{W}\beta_2,
\end{align}
where $\beta_2$ is the last $(p-q)$ elements of $\beta$.

Luckily, when the IID assumption holds, this bias can often be mitigated by a shift in the residual terms used to construct the predictive interval. Note that,  the expectations of the residual terms are not zero: 
\begin{align}
    \mathbb{E} \left\{ (1-g_{i,new})v_i |\mathbf{X}, \X_{new}   \right\}& = \frac{1-g_{i,new}}{1-g_{ii}}\left( \W_i^T-\Z_i^T(\mathbf{Z}^T\mathbf{Z})^{-1}\mathbf{Z}^T\mathbf{W}\right)\beta_2 \nonumber\\
    & = \frac{1-g_{i,new}}{1-g_{ii}} (\W_i^\perp)^T \beta_2 \overset{def}{=} shift
\end{align}
where $\W_i^\perp$ is the $i$th row of the matrix $\mathbf{W}^\perp = \{I - \mathbf{Z}(\mathbf{Z}^T\mathbf{Z})^{-1} \mathbf{Z}^T\} \mathbf{W}$.
The $shift$ and $bias$ often have the opposite signs and thus, when added together,  
they cancel each other to a certain extent.

For example, suppose a new individual case is an ``average individual" of the training data with $\X_{new} = \bar \X = \frac{1}{n}\sum_{i = 1}^n \X_i$. Then, the bias of the point predictor and  the average shift of the residual terms are 
$
bias = - (\bar\W^\perp)^T \beta_2$ and $average \,\, shift = \frac{1}{n}\sum_{i = 1}^n \left\{ \frac{1-{1}/{n}}{1-g_{ii}} (\W_i^\perp)^T \beta_2\right\},  
$
respectively. 
Since $\frac{1 - {1}/{n}}{1 - g_{ii}} \approx 1$ when $\Z_i$'s are IID (cf., Lemma A1 in Supplementary), the {\it average shift} $\approx (\bar\W^\perp)^T \beta_2 = -${\it bias}, thus they are approximately canceled out 
in the predictive interval~(\ref{eq:PI-wrong}). 
This cancellation explains in part why the prediction interval (\ref{eq:PI-wrong}) is still roughly on target, even if the learning model is wrong. The cancellation is not as complete, when the testing data $\X_{new}$ is just an IID sample and not the ``average'' $\bar \X$. 
It appears that the combination of an enlarged interval and the cancellation of the {\it bias} and {\it shift} helps  ensure the validity of conformal prediction under a wrong model for IID testing data. This self balance to maintain validity mirrors a homeostasis process and we referred to it as a {\it homeostasis phenomenon}. 

A wrong learning model also has implications on the lengths of the prediction intervals. 
The proposition below states that the width of the predictive interval based on the wrong model $\mu_1(\cdot)$ is expected to be wider than that based on the correct model $\mu_0(\cdot)$, if $\X_{new} = \bar \X$. 
A proof can be found in the Supplementary. 

\begin{prop} 
\label{Prop:1} Under model (\ref{eq:LM-mu0}), assume $\epsilon_i \overset{iid}{\sim}N(0, \sigma^2)$, $\X_i$'s are IID from a normal distribution and $\beta_2^T \Sigma_{\rm w|z} \beta_2 > 0$, where $\Sigma_{\rm w|z} = {\rm var}(\W_1| \Z_1)$. 
Suppose $\X_{new}=\bar \X$, then
    \begin{align*}
    &\lim_{n\to\infty} \mathbb{P}[ q_{1-\frac{\alpha}{2}}\left(\{(1-g_{i,new})v_i\}_{i=1}^n  \right)- q_{\frac{\alpha}{2}}\left(\{(1-g_{i,new})v_i\}_{i=1}^n \right)\\
    &\quad\quad\quad > q_{1-\frac{\alpha}{2}}\left(\{(1-h_{i,new})u_i\}_{i=1}^n  \right)- q_{\frac{\alpha}{2}}\left(\{(1-h_{i,new})u_i\}_{i=1}^n \right) \big] =1.
    \end{align*}
    That is, with probability tending to 1, the width of predictive interval (\ref{eq:PI-true}) $\geq$ the width of predictive interval (\ref{eq:PI-wrong}).  
\end{prop}

The following  numerical example provides empirical evidence to support our discussions. 

\begin{example} Suppose 
we have only two covariates in the true model 
\begin{align}
\label{eq:m0}
    y_i = \mu_0(\X_i) +  \epsilon_i =  \beta_0+\beta_1 z_{i}+\beta_2 w_{i} +\epsilon_i, \quad \epsilon_i \overset{\tiny iid}{\sim} N(0,\sigma^2)
\end{align}
where $\X_i = (z_{i}, w_{i})^T \overset{\tiny iid}{\sim} N(\mu_x,\Sigma_x)$ and $\epsilon_i$ and $\X_i$ are independent. In our numerical study, 
$(\beta_0, \beta_1, \beta_2) = (-1,2,2)$, $\sigma^2 =1$, $\mu_x = (0,0)^T$, the $(k, k')$-element of $\Sigma_x$ is 
$0.5^{|k-k'|}/2$, $k, k' \in \{1,2\}$ and $n = 300$. 

For the testing data, we consider two scenarios: (i) under the IID assumption that ${\X}_{new} \sim N(\mu_x,\Sigma_x)$ and $({\X}_{new}, y_{new})$ follows (\ref{eq:m0}); (ii) the marginal distribution of ${\X}_{new}$ is instead from ${\X}_{new} \sim N(\tilde\mu_x,\tilde\Sigma_x)$ and, given ${\X}_{new}$, $({\X}_{new}, y_{new})$ follows (\ref{eq:m0}). Here, $\tilde\mu_x = \mu_x + (2, 2)^T$ and 
the $(k, k')$-element of $\tilde\Sigma_x$ is 
$0.8^{|k-k'|}/2$, $k, k' \in \{1,2\}$.

In addition to the correct model (a) $\mu_0(\X_i) = \beta_0+\beta_1 z_{i}+\beta_2 w_{i}$, three wrong learning models are considered:
\begin{align*}
\hbox{(b)} \quad  \mu_1(\X_i) & = \gamma_0+\gamma_1 z_{i} \quad \hbox{(partially correct, without covariate $w_i$);} \\
 \hbox{(c)} \quad \mu_2(\X_i) & = \xi_0+\xi_1 z_{i}^2 \quad \hbox{(a wrong regression form).} \\
 \hbox{(d)} \quad \mu_3(\X_i) & = \eta_0 \quad \hbox{ (without any covariates);}
\end{align*}
For model fitting, we use the least squares method in all three cases. 

Reported in each cell of Table 1 are the coverage rate and average length (inside brackets) of $95\%$  conformal predictive intervals for $y_{new}$, computed based on $200$ repetitions. As expected, in the IID scenario, all learning models can provide valid prediction results. However, the smallest interval length   is observed under the true model.
In the non-IID scenario, only the true model can provide a valid prediction. The other three learning models do not provide valid predictive inference in terms of a correct  coverage rate, even though their predictive intervals are wider. 
The results in both scenarios underscore the importance of using a correct learning model for prediction.

\begin{table}
{\scriptsize
\begin{center} \centering
\begin{tabular}{c|c|c|c|c}
\hline\hline
&\multicolumn{1}{|c|}{True model}&\multicolumn{3}{|c}{Wrong model}
\\ \hline
 & $\mu_0(\cdot)$ & $\mu_1(\cdot)$ & $\mu_2(\cdot)$ & $\mu_3(\cdot)$ \\
\hline  
IID Scenario &.985 (4.420) & .96 (6.957) &.98 (11.697) &.98 (12.147)\\
\hline 
Non-IID Scenario &.985 (4.421) &.81 (6.957) &.345 (11.280) &.33 (12.147) \\
\hline 
\end{tabular}
\end{center} }
\caption{\small Performance of $95\%$ predictive intervals under four learning models and in two scenarios (coverage rates (before brackets) and average interval lengths (inside brackets)). Model $\mu_1(\cdot)$ is a partially wrong model, $\mu_2(\cdot)$ is a completely wrong model and $\mu_3(\cdot)$
does not use any covariates.  Training data size = 300; Testing data size = 1; Repetition = 200.}
\end{table}

In order to get the full picture of the predictive intervals of all confidence levels 
under different scenarios and different learning models, we plot in Figure~\ref{fig:PV_LM} the predictive curves obtained from the first realization of the $200$ repetitions (Other $199$ realizations produce more or less the same plots). 
Plot (a) is for $\X_{new} = (-0.011,-0.046) = \frac{1}n \sum_{i=1}^n \X_i$, (b) is for $\X_{new}=(0.365,-0.026) \overset{iid}{\sim} \X_i$ and (c) $\X_{new}=(2.44,2.09) \not \sim \X_i$. In each plot, we have four predictive curves corresponding to four working models, plus the target (oracle) predictive curve of $PV(y) = 2 \max\{\Phi(y - \mu_{new}), 1 - \Phi(y - \mu_{new})\}$  obtained by pretending that we know exactly $y_{new}$'s distribution: $y_{new} \sim N(\mu_{new}, 0.5)$ with $\mu_{new} = (-1,2,2)\X_{new}$.
In each of the plots (a)-(c), the predictive curves trained with the correct model (black solid curves) are very close to the target oracle predictive curves (red solid curves), indicating that if we use the true model as the learning model, we are able to provide very accurate prediction at all  confidence levels.
Under the wrong models, however, the take-home messages are very different. 
In plot (a) with $\X_{new}$ being the ``average individual,'' we see an almost complete cancellation of $bias$ and $shift$ as described earlier. However, the predictive curves are much wider than those based on the correct model.
Plot (b) is for the IID case of $\X_{new} \sim \X_i$. In this case the curves are similar to those in plot (a), although the cancellations are not as complete as for the `average individual'. Nevertheless, the enlarged interval widths help maintain the coverage.
Plot (c) is for non-IID case, in which the cancellations of $bias$ and $shift$ are not effective when wrong learning models are used, leading to wrong predictions. 
In plots (a) - (c), we can also see that a partially correct model $\mu_1(\cdot)$ performs better than the other two completely wrong models $\mu_2(\cdot)$ and $\mu_3(\cdot)$. 

In summary, when we train prediction algorithms  using a wrong model, the IID assumption is essential for the validity of prediction, and using a wrong model often results in wider, sometimes much wider, predictive intervals. When we train the same algorithms using the correct model, the validity and efficiency of the predictions are observed in both IID and Non-IID scenarios conditional on $\X_{new}$.

\begin{figure}
    \centering
\includegraphics[width=3.9cm,height=4.5cm]{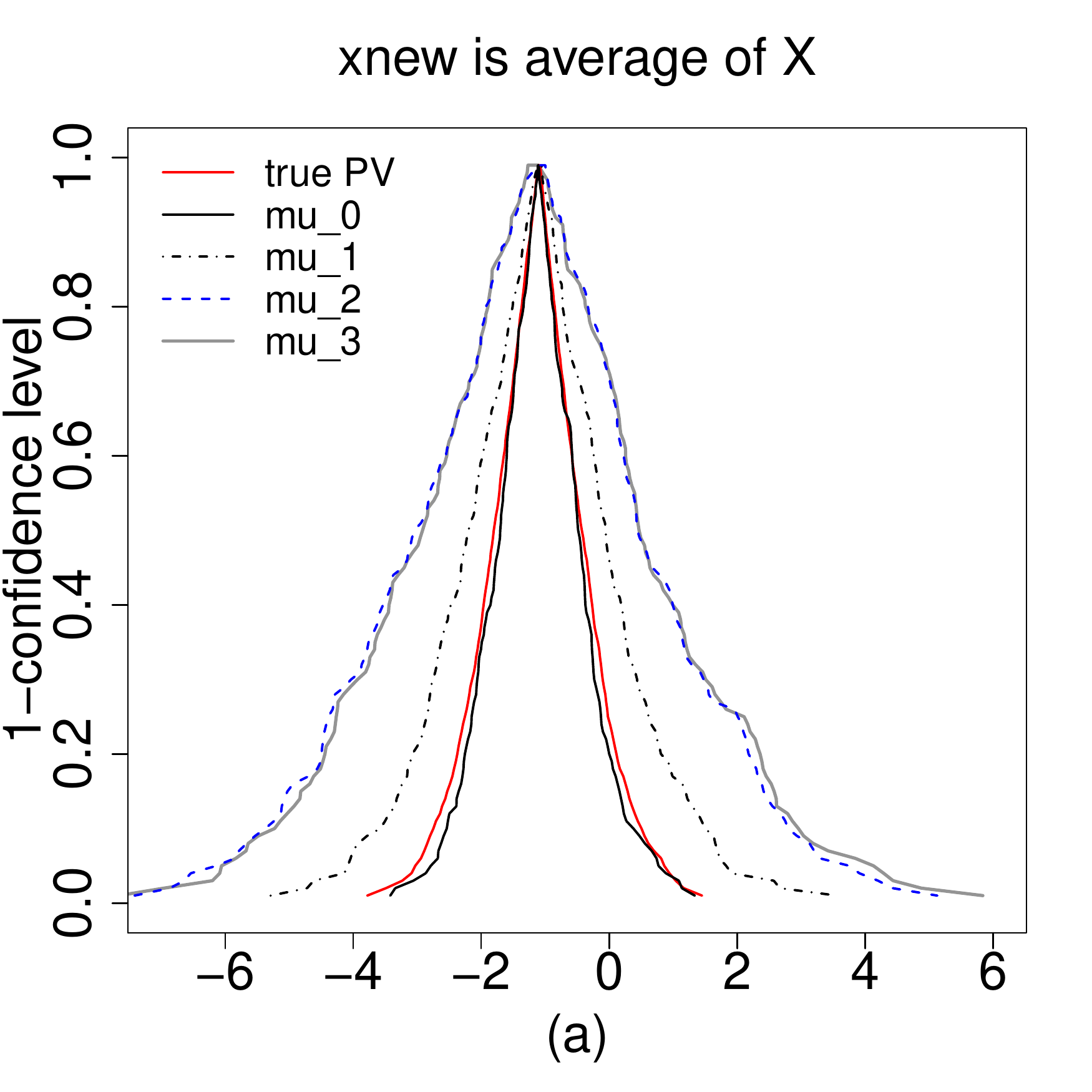}
\includegraphics[width=3.9cm,height=4.5cm]{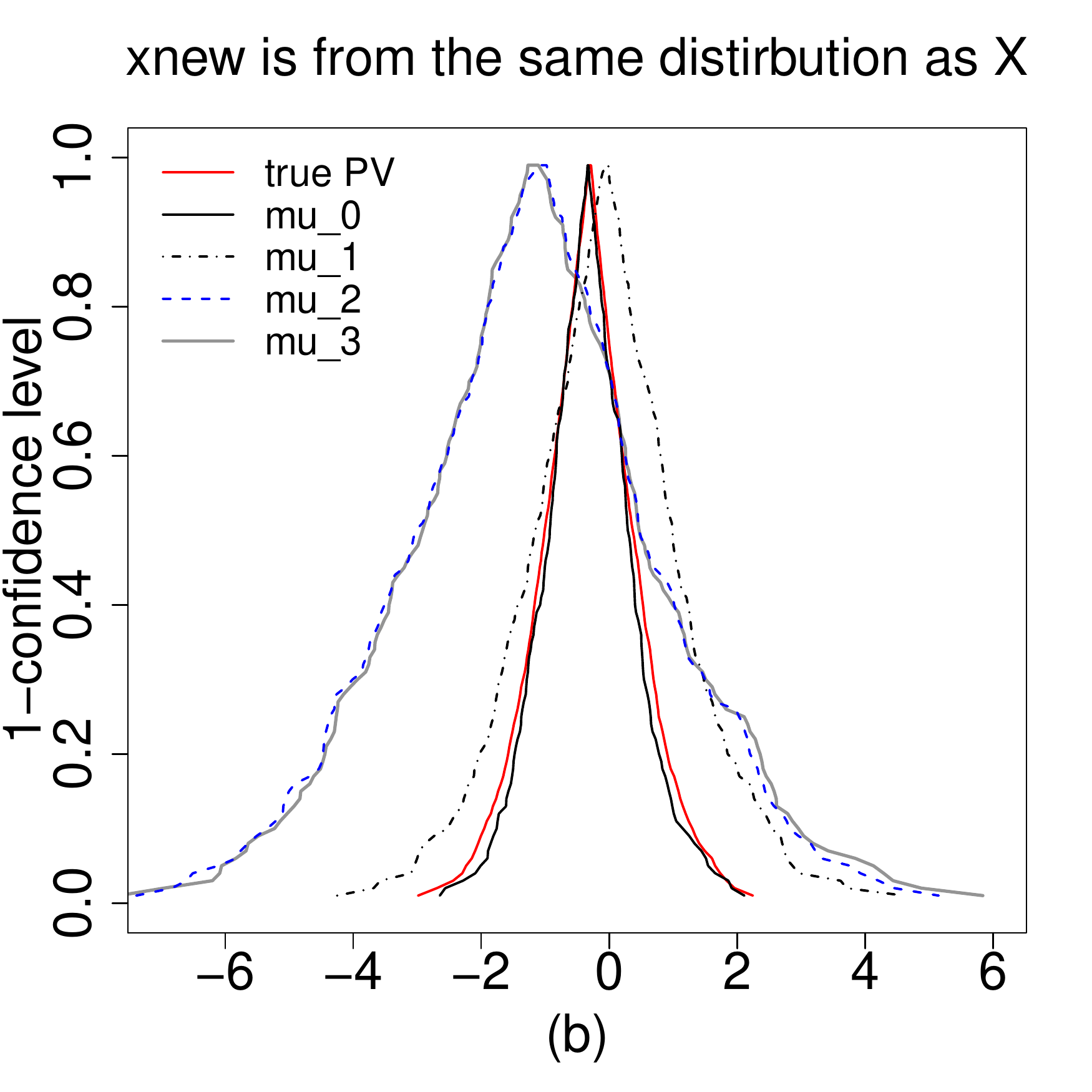}
\includegraphics[width=3.9cm,height=4.5cm]{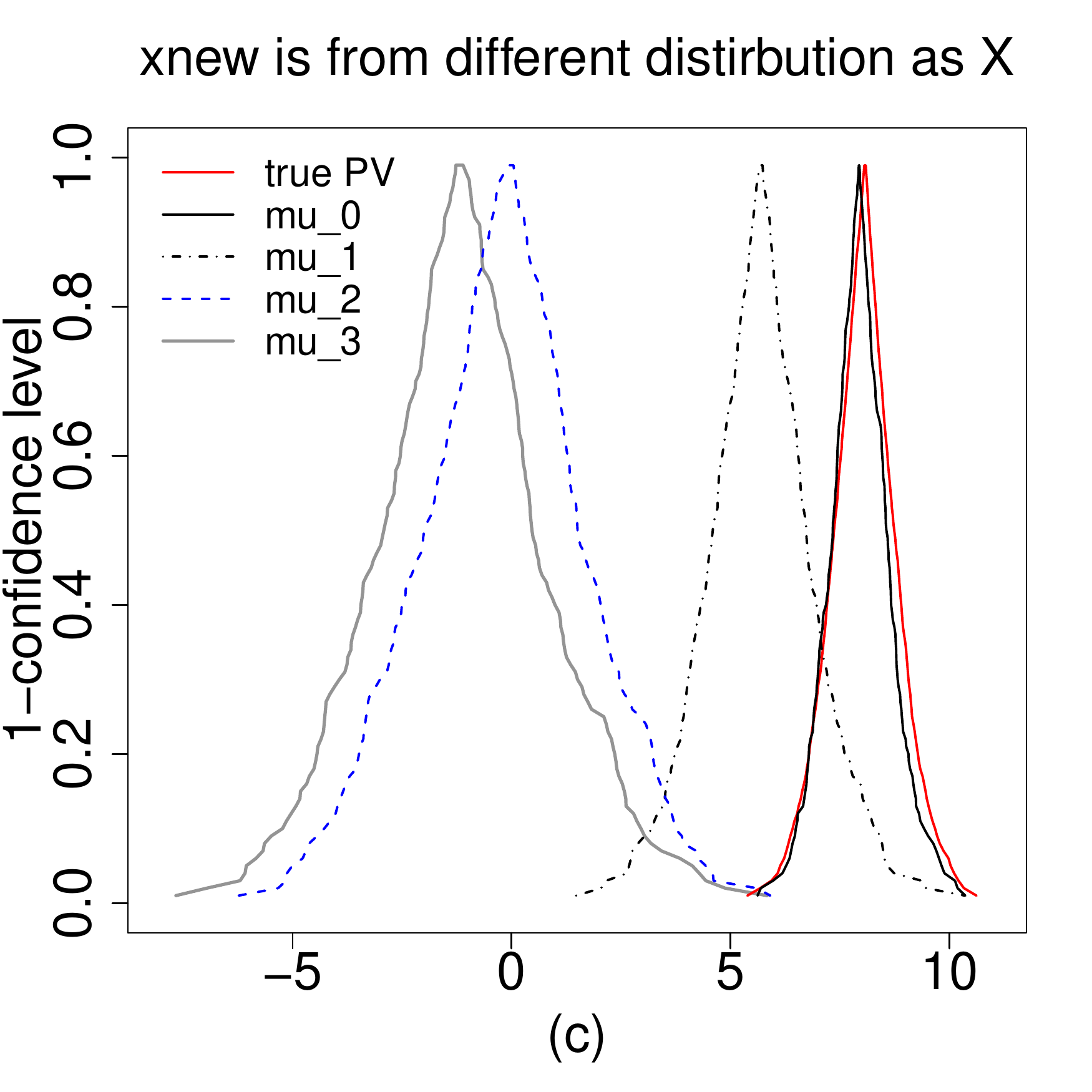}
 \caption{\small Plots of predictive curves for (a) $\X_{new} = \bar \X = \frac{1}n \sum_{i=1}^n \X_i$; (b)  $\X_{new} \overset{iid}{\sim} \X_i$ and (c) $\X_{new} \not \sim \X_i$. In each plot, the red solid curve is the target (oracle) predictive curve  $PV_n(y) = 2 \max\{\Phi(y - \mu_{new}), 1 - \Phi(y - \mu_{new})\}$, obtained assuming that the distribution of $y_{new}$ is completely known.
 The predictive curves in black and blue are obtained using the four working models, respectively. The solid black curve is for learning model $\mu_0(\cdot)$, the dotted black for $\mu_1(
 \cdot)$, dashed blue for  $\mu_2(\cdot)$ and solid gray for $\mu_3(\cdot)$.
 }
  \label{fig:PV_LM}
\end{figure}

\end{example}

\subsection{Prediction in neural network model} 

The discussion in the linear model in Section~\ref{sec:3.1} can be extended to other models. We consider in this subsection an example of simple neural network models. We use a simulation study to provide empirical support for our discussion. Note that, in the current neural network development, model fitting algorithms do not pay much attention to correctly estimate the model parameters. In addition to what we learned in the linear model, we find that the estimation of  model parameters plays an important role in prediction as well.

\begin{example}
Suppose our training data $(y_i, \X_i)$, $i = 1, \ldots, n$, are IID samples from the model
\begin{equation}
\label{eq:NN0}
     y_i = \mu_0(\X_i) + \epsilon_i= \max\big\{0,\max\{0, \, z_{i1}+ z_{i2}\}-\max\{0,w_i\}\big\} + \epsilon_i, \quad \epsilon_i \overset{\tiny iid}{\sim} N(0, \sigma^2), 
\end{equation}
where $\X_i = (z_{i1}, z_{i2}, w_{i})^T \overset{\tiny iid}{\sim} N(\mu_x,\Sigma_x)$ and $\epsilon_i$ and $\X_i$ are independent. Here, $\mu_x = (0,0, 0)^T$, the $(k, k')$-element of $\Sigma_x$ is 
$0.5^{|k-k'|}/2$, for $k, k' \in \{1,2, 3\}$, $\sigma^2 =1$ and $n = 300$. 
Model (\ref{eq:NN0}) is in fact a neural network model (with a diagram presented in Figure~\ref{fig:NN-d} (a)) and we can re-express~$\mu_0(\X_i)$~as
\begin{equation}
  \mu_0(\X_i) = f\big(A_2 f(A_1\X_i)\big)
\label{eq:NN-mu0}  
\end{equation}
Here, $f(x) = \max(x,0)$ is the ReLU activate function, and 
$A_1 = \left(\begin{array}{ccc}
         a_{11}^{(1)}& a_{12}^{(1)}&a_{13}^{(1)} \\
         a_{21}^{(1)}& a_{22}^{(1)}&a_{23}^{(1)} 
    \end{array}\right)$ and $A_2 = \big( a_1^{(2)},
         a_2^{(2)} \big)$ are the model parameters. 
Corresponding to (\ref{eq:NN0}), the true model parameter values are $a_{11}^{(1)} = a_{12}^{(1)} = a_{23}^{(1)} = 1$, $a_{13}^{(1)} = a_{21}^{(1)} = a_{22}^{(1)} = 0$ and $(a_1^{(2)}, a_2^{(2)}) = (1,-1)$. 
In our analysis, we assume that we know the model form (\ref{eq:NN-mu0}) but do not know the values of model parameters $A_1$ and $A_2$.

\begin{figure}
    \centering
    \includegraphics[height=4.8cm,width=14cm]{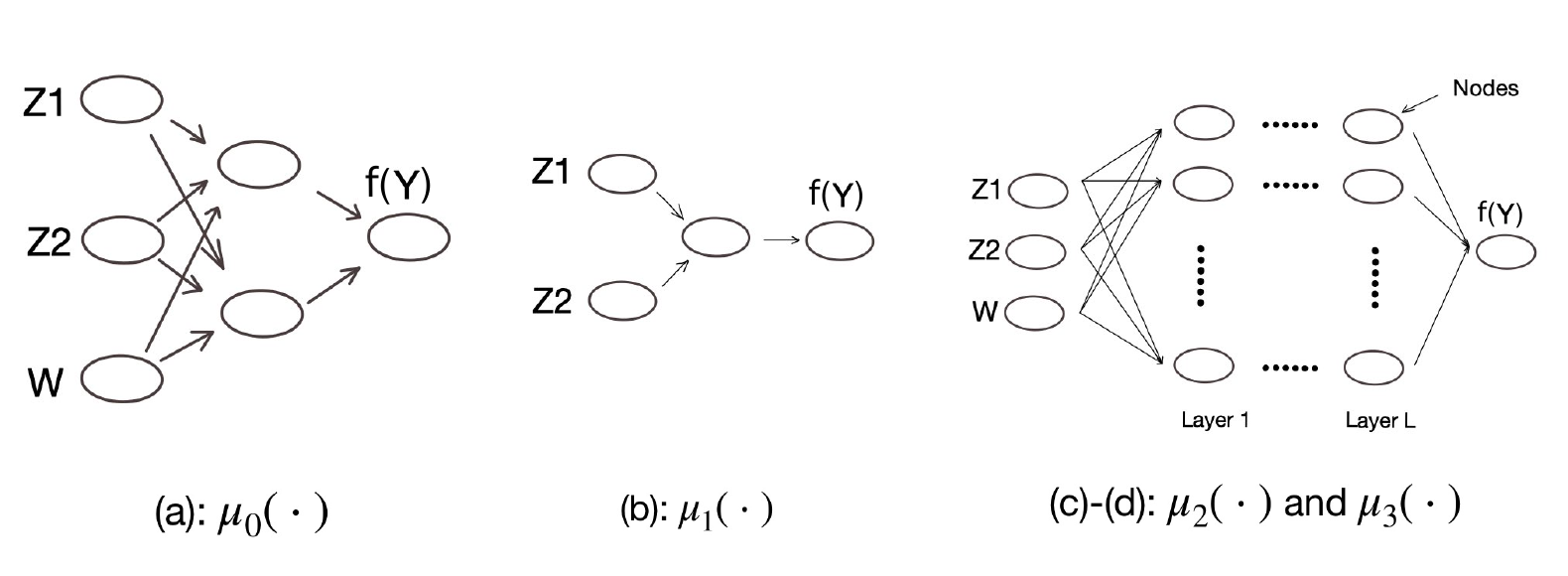}
    \vspace{-6mm} 
    \caption{\small Diagrams of four neural network models: (a) true $\mu_0(\cdot)$; (b) partial $\mu_1(\cdot)$; and (c)-(d) over-parametrized  $\mu_2(\cdot)$  and $\mu_3(\cdot)$ of ($20$~nodes in each~layer) $\times$ $L$ layers, with $L= 20$ and $100$, respectively.}
    \label{fig:NN-d}
\end{figure}

For the testing data, we consider two scenarios: (i) [IID case] ${\X}_{new} \overset{\tiny iid}{\sim} N(\mu_x,\Sigma_x)$ and, given ${\X}_{new}$, $({\X}_{new}, y_{new})$ follows (\ref{eq:NN0}); (ii) [Non-IID case] the marginal distribution 
${\X}_{new} \overset{\tiny i.i.d}{\sim} N(T_1,T_2,T_3)$ and, given ${\X}_{new}$, $({\X}_{new}, y_{new})$ follows (\ref{eq:NN0}). Here, $T_1,T_2,T_3$ are i.i.d random variables from t distribution with degrees of freedom 3 and non-centrality parameter 1. 

In addition to (a) the true model $\mu_0(\cdot)$,
four wrong learning models are considered:
\begin{align} 
 & \hbox{(b)} \quad \mu_1(\X_i) = f(B\Z_i)
 \quad \hbox{\rm (partially correct neural network model, missing $w_i$);} \nonumber 
 \\
    &\hbox{(c)} \quad  \mu_2(\X_i) = f\big(C_{20}f(C_{19}\cdots f(C_1\X))\big) \qquad \hbox{\rm (deep neural network model with 20 layers)}; 
    \nonumber \\
    &\hbox{(d)} \quad  \mu_3(\X_i) = f\big(D_{100}f(D_{99}\cdots f(D_1\X))\big) \qquad \hbox{\rm (deep neural network model with 100 layers)}; \nonumber
    \\
    &\hbox{(e)} \quad  \mu_4(\X_i) = \eta_0 \qquad \hbox{\rm (without any covariates)}, \nonumber 
\end{align}
where $\Z_i = (z_{i1},z_{i2})^T$, $B = (b_{1}, b_{2})$, $C_1,D_1 \in \mathbb{R}^{20\times 3}$, $C_{20},D_{100} \in \mathbb{R}^{1 \times 20}$, and $C_i,D_j\in \mathbb{R}^{20\times 20}$, $2\leq i \leq 19$, $2\leq j \leq 99$.
In our analysis, the neural network models $\mu_0(\cdot)$ - $\mu_3(\cdot)$ are fitted using the {\sc neuralnet} package (\hyperref[https://cran.r-project.org/web/packages/neuralnet/]{cran.r-project.org/web/packages/neuralnet/}).

The Neuralnet package is an off-the-shelf machine learning algorithm. 
Its emphasis is on learning and not on model parameter estimation. Even under the true model $\mu_0(\cdot)$, the estimates of model parameters from Neuralnet are not very accurate; cf., Table~\ref{tab:MSE}.
In the table, ``Opt-MSE'' refers to a code that we wrote by directly minimizing {\it MSE} $= \sum_{j=1}^n (y_j - \mu_0(\X_j))^2$, which can be implemented when the neural network is small. The calculation is based on $20$ repeated runs, each with a training data set of size $n = 300$ from model~(\ref{eq:NN0}). 

\begin{table}[]
    \centering
    {\scriptsize
    \begin{tabular}{c|c|c|c|c|c|c|c|c}
        \hline\hline
          MSE&$a_{11}$&$a_{12}$&$a_{13}$&$a_{21}$&$a_{22}$&$a_{23}$&$b_1$&$b_2$  \\
        \hline
          Opt-MSE& 0.07&0.059&0.31&0.154&0.109&0.124&0.06&0.101\\
         \hline
         Neuralnet&4.87&5.9&1.53&1.14&2.23&2.08&4.87&0.84\\
         \hline
    \end{tabular}
 }   \caption{\small
      Mean square error of each parameter in $\mu_0$ using two estimation procedures (Training data $n = 300$; Repetition = $10$)
     }
    \label{tab:MSE}
\end{table}

Reported in Table~\ref{tab:NN} are the coverage rate and average interval length of predictive intervals computed under $10 = 5 \times 2$ settings with five different learning models $\mu_k(\cdot)$, $k = 0,1,\cdots, 4$, and in two scenarios.
The analysis is repeated for $10$ times with $10$ simulated training datasets from model (\ref{eq:NN0}). We use $10$ repetitions and not a greater number, because it takes a long time to fit a neural network model. However, for each of the $10$ training data sets,  $20$ pairs of $(y_{new}, \X_{new})$ are used. So, for the reported values, 
each is computed using $10 \times 20 = 200$ pairs of $(y_{new}, \X_{new})$. 
For the true neural network model $\mu_0(\cdot)$, Opt-MSE is also used to fit the model. 
As we can see in Table~\ref{tab:NN}, under the IID Scenario, all predictive intervals are valid with a  correct coverage. The best one with the shortest interval length is the one that uses the correct model and Opt-MSE estimation method. In the non-IID case, only the shallow neural network models provide valid predictions, and among them, Opt-MSE methods can give us confidence intervals with half the width.
Indeed, when a wrong learning model is used, the IID assumption is essential for the prediction validity and the use of a wrong model often results in wider intervals. Furthermore, the estimation of model parameters seems to also have big impact on prediction.  

\begin{table}[]
    \centering
{\scriptsize
\begin{tabular}{c|c|c|c|c|c|c}
\hline\hline
 &\multicolumn{2}{|c|}{True model}&\multicolumn{4}{|c}{Wrong model} \\
\hline
 &\multicolumn{2}{|c|}{$\mu_0(\cdot)$} &{ $\mu_1(\cdot)$}&{$\mu_2(\cdot)$}&{$\mu_3(\cdot)$}&{$\mu_4(\cdot)$}\\
\hline
 &Opt-MSE &Neuralnet &Nueralnet &Nueralnet &Nueralnet &Nueralnet
 \\
\hline  
IID Scenario &.995 (4.462) &.99 (4.608)& .99 (4.809) & .99 (5.212)& .99 (5.201) &.985 (5.26)\\
\hline 
Non-IID Scenario &.955 (4.52) &.985 (9.327) &.98 (9.77) & .71 (5.899) & .695 (5.201) & .685 (5.277)\\
\hline 
\end{tabular}
}
    \caption{ \small Performance of $95\%$ predictive intervals under five different learning models and in two scenarios: coverage rates (before brackets) and average interval lengths (inside brackets) (Training data size = 300; Testing data size = 20; Repetition = 10)}
\label{tab:NN}
\end{table}

To get a full picture of the predictive intervals at all levels,
we plot in Figure~\ref{fig:PV-NN} the predictive curves of $y_{new}$.
The plots  are based on the first training dataset and making prediction for (a) the IID case with the realizatiof ${\bf x}_{new} = (-0.909,-1.149,-0.771)$, 
and (b) the non-IID case with the realization ${\bf x}_{new} = (3.653,1.748,1.063)$.
The realized value of  $\mu_0({\bf x}_{new})$ is $0$ and $4.338$ in (a) and (b),  respectively. 
From Figure~\ref{fig:PV-NN}, we see 
that the use of a wrong model $\mu_1(\cdot)$ - $\mu_4(\cdot)$ 
results in wider predictive curve
(and predictive intervals at all levels $1 -\alpha \in (0,1)$) in both IID and non-IID cases.
Although the shallow neural network models 
$\mu_0(\cdot)$ and $\mu_1(\cdot)$ can provide good coverage rates, the predictive curves in the non-IID case are much fatter than other approaches. This peculiar phenomenon occurs even when we assume to know the true model structure $\mu_0(\cdot)$, indicating the importance of  estimating model parameters accurately. 
Furthermore, in the non-IID case, there are large shifts
when using  deep neural network models
$\mu_2(\cdot)$ and $\mu_3(\cdot)$, 
leading to invalid predictions. 
The best prediction result is from
the one obtained by using the correct learning model $\mu_0(\cdot)$ with the more accurate parameter estimation method Opt-MSE.
The message is the same as what we have learned from Table~\ref{tab:NN}. In addition to what we learned in the linear setting in Section~\ref{sec:3.1}, 
a good estimation of model parameters is important in prediction.

\begin{figure}
    \centering
\includegraphics[width=8.4cm,height=5.2cm]{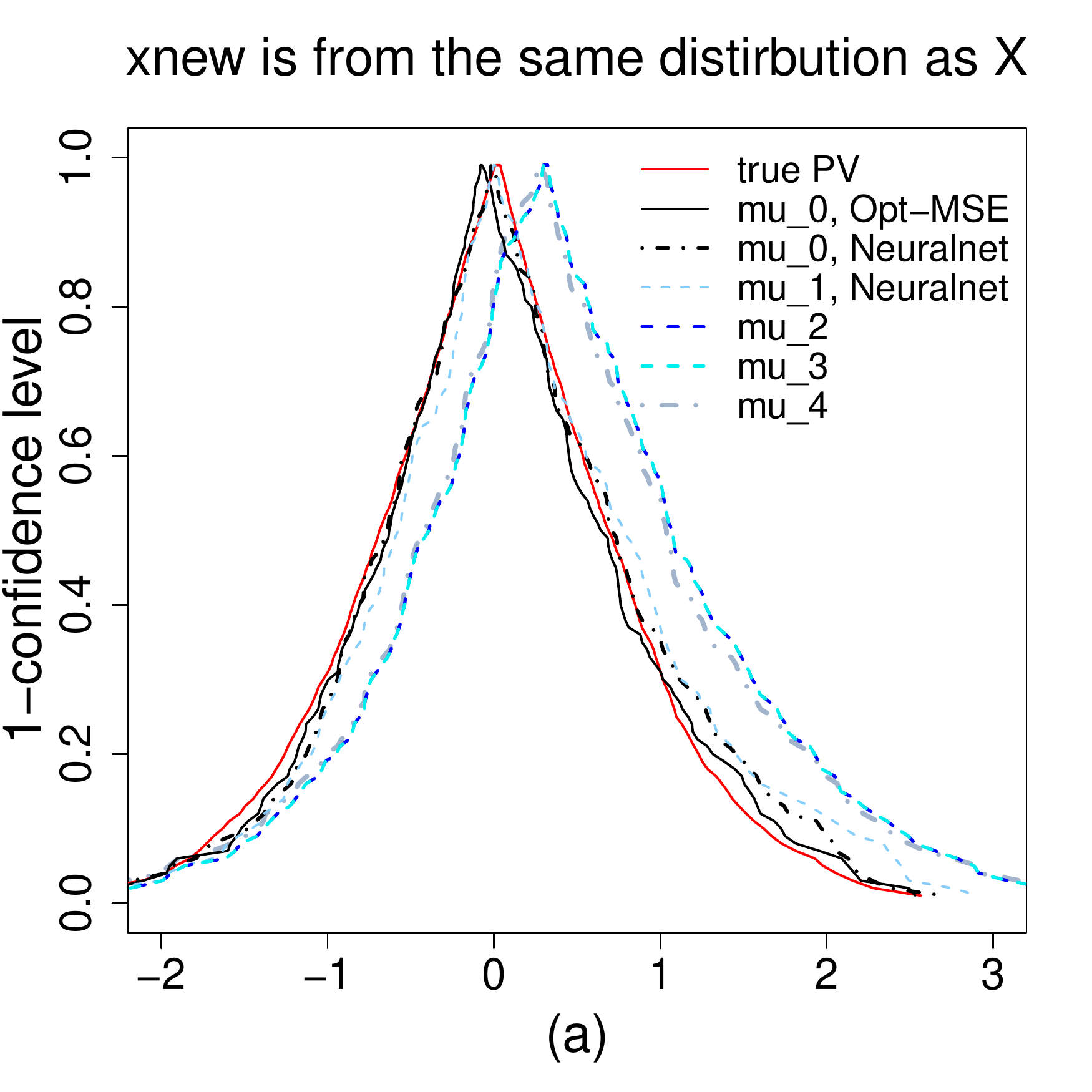}
\includegraphics[width=8.4cm,height=5.2cm]{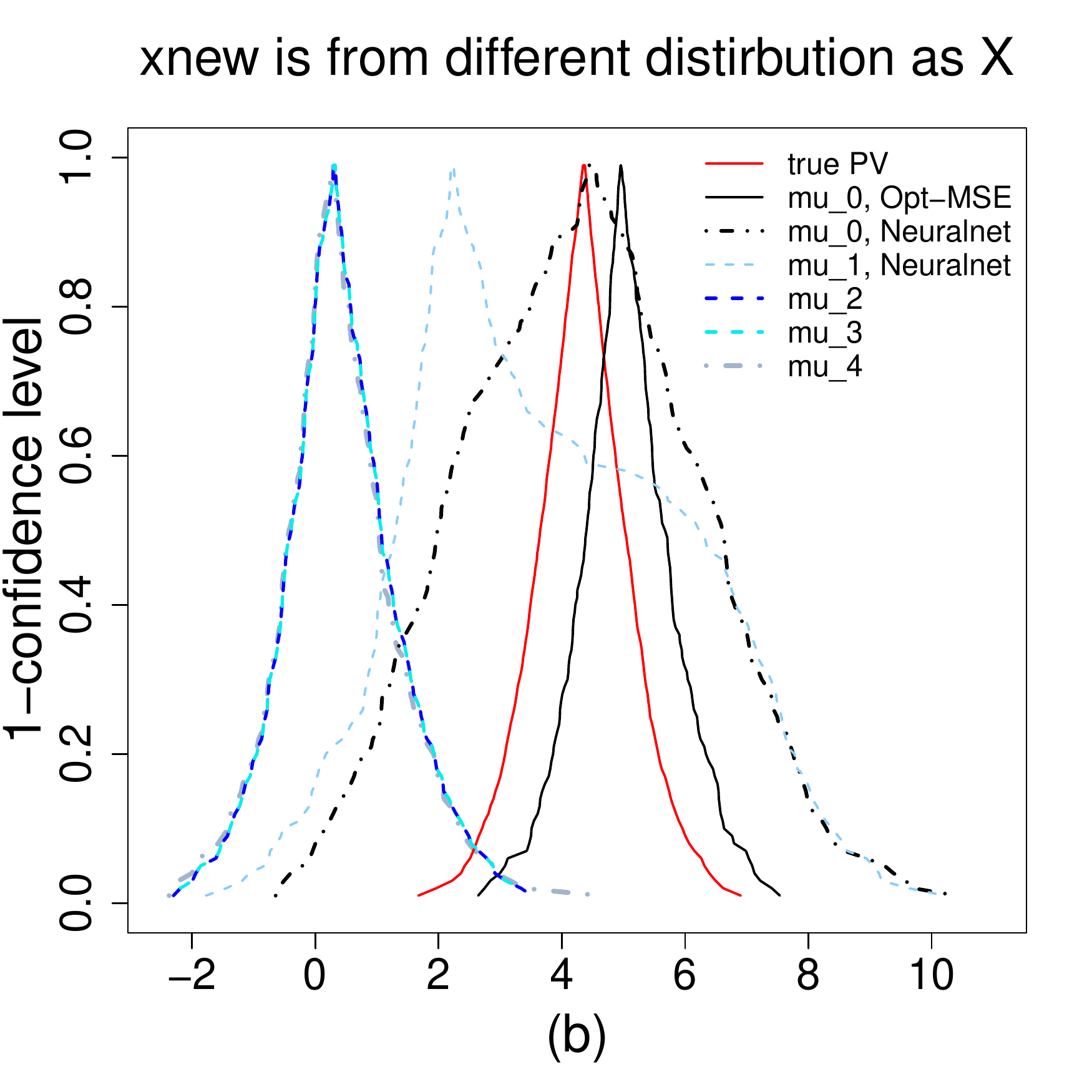}
\caption{\small 
Plots of predictive curves for (a)  $\X_{new} \overset{iid}{\sim} \X_i$ and (b) $\X_{new} \not \sim \X_i$. In each plot, the red solid curve is the target (oracle) predictive curve  $PV_n(y) = 2 \max\{\Phi(y - \mu_{new}), 1 - \Phi(y - \mu_{new})\}$, obtained assuming that the distribution of $y_{new}\sim N(\mu_{new}, 1)$ is completely known.
 The two predictive curves obtained using $\mu_0(\cdot)$ are in black  (solid line for Opt-MSE; dashed line for Neuralnet). The other predictive curves (all in a dashed or broken line and in various colors) are obtained using the other four wrong working models.
 }
    \label{fig:PV-NN}
\end{figure}
\end{example}

\section{Conclusion}
Professor Efron pointed out that
``the 21st Century has seen the rise of a new breed of what can be called `pure prediction
algorithms.'" We are fully in agreement with Professor Efron's
discussion that the
prediction algorithms ``can be stunningly successful,'' and that 
``the emperor
has nice clothes but they're not suitable for every occasion.'' 
Along the same line and under the setting of conformal prediction, we have demonstrated and explained  how and why a prediction method can be successful under the IID assumption, even if the learning model is completely 
wrong. More importantly, we have also demonstrated that it is still
meaningful, and often crucial, to build our prediction algorithms based on a good practice of
modeling, estimation and inference.

Learning is a powerful tool for processing large data for information in modern data science. The impressive narrative of the black-box approaches thus far, however, is only a step in our long journey ahead for the Twenty-first Century
statistics and data science. We have provided an explanation why current machine learning methods work so well in academic research settings (i.e., randomly splitting the training and testing data) but fail to deliver as successful stories in applications.  We believe that lacking ingredients are ``the most powerful ideas of Twentieth Century
statistics'' ---
modeling, estimation and inference. We are excited about our profession. 
We fully anticipate
and believe that statistical developments, with rich traditions, will play a pivotal role in building the mathematical foundation of modern data science and in fully realizing its potential for real-world applications.

\bibliographystyle{asa}
\bibliography{reference}

\begin{thebibliography}{14}
\newcommand{\enquote}[1]{``#1''}
\expandafter\ifx\csname natexlab\endcsname\relax\def\natexlab#1{#1}\fi

\bibitem[{Barber et~al.(2019{\natexlab{a}})Barber, Candes, Ramdas, and
  Tibshirani}]{barber2019limits}
Barber, R.~F., Candes, E.~J., Ramdas, A., and Tibshirani, R.~J.
  (2019{\natexlab{a}}), \enquote{The limits of distribution-free conditional
  predictive inference,} \textit{arXiv preprint arXiv:1903.04684}.

\bibitem[{Barber et~al.(2019{\natexlab{b}})Barber, Candes, Ramdas, and
  Tibshirani}]{barber2019predictive}
--- (2019{\natexlab{b}}), \enquote{Predictive inference with the jackknife+,}
  \textit{arXiv preprint arXiv:1905.02928}.

\bibitem[{Birnbaum(1961)}]{Birnbaum1961}
Birnbaum, A. (1961), \enquote{Confidence curves: An omnibus technique for
  estimation and testing statistical hypotheses,} \textit{Journal of the
  American Statistical Association}, 56, 246--249.

\bibitem[{Cox(1958)}]{cox1958}
Cox, D.~R. (1958), \enquote{Some problems connected with statistical
  inference,} \textit{The Annals of Mathematical Statistics}, 29, 357--372.

\bibitem[{Efron(1993)}]{Efron1993}
Efron, B. (1993), \enquote{Bayes and likelihood calculations from confidence
  intervals,} \textit{Biometrika}, 80, 3–26.

\bibitem[{Efron(2020)}]{Efron2020}
--- (2020), \enquote{Prediction, Estimation, and Attribution (with
  discussion),} \textit{Journal of the American Statistical Association}, To
  appear.

\bibitem[{Lawless and Fredette(2005)}]{Lawless2005}
Lawless, F. and Fredette, M. (2005), \enquote{Frequentist prediction intervals
  and predictive distributions,} \textit{Biometrika}, 92, 529--542.

\bibitem[{Lei et~al.(2018)Lei, G’Sell, Rinaldo, Tibshirani, and
  Wasserman}]{lei2018distribution}
Lei, J., G’Sell, M., Rinaldo, A., Tibshirani, R.~J., and Wasserman, L.
  (2018), \enquote{Distribution-free predictive inference for regression,}
  \textit{Journal of the American Statistical Association}, 113, 1094--1111.

\bibitem[{Schweder and Hjort(2016)}]{Schweder2016}
Schweder, T. and Hjort, N. (2016), \textit{Confidence, Likelihood and
  Probability}, Cambridge, U.K.: Cambridge University Press.

\bibitem[{Shafer and Vovk(2008)}]{shafer2008tutorial}
Shafer, G. and Vovk, V. (2008), \enquote{A tutorial on conformal prediction,}
  \textit{Journal of Machine Learning Research}, 9, 371--421.

\bibitem[{Shen et~al.(2018)Shen, Liu, and Xie}]{Shen2018}
Shen, J., Liu, R., and Xie, M. (2018), \enquote{Prediction with confidence—A
  general framework for predictive inference,} \textit{Journal of Statistical
  Planning and Inference}, 195, 126--140.

\bibitem[{Vovk et~al.(2005)Vovk, Gammerman, and Shafer}]{vovk2005algorithmic}
Vovk, V., Gammerman, A., and Shafer, G. (2005), \textit{Algorithmic learning in
  a random world}, Springer Science \& Business Media.

\bibitem[{Vovk et~al.(2019)Vovk, Shen, Manokhin, and Xie}]{Vovk2019}
Vovk, V., Shen, J., Manokhin, V., and Xie, M. (2019), \enquote{Nonparametric
  predictive distributions by conformal prediction,} \textit{Machine Learning},
  108, 445--474.

\bibitem[{Xie and Singh(2013)}]{Xie2013}
Xie, M. and Singh, K. (2013), \enquote{Confidence distribution, the frequentist
  distribution estimator of a parameter (with discussion),}
  \textit{International Statistical Review}, 81, 3--39.

\end{thebibliography}

\section*{Supplementary: Proof of Propositions 1 and 2}

\subsection*{S.1 \, Proof of Proposition 1.}
{\bf Proof of Proposition 1.}
We only prove that 
$\mathbb{P}(y_{new}\in \{y: Q_n(y) \ge \frac{\alpha}{2}\}) \geq 1 - \alpha$. The proof of $\mathbb{P}(y_{new}\in \{y: 1 - Q_n(y) \ge \frac{\alpha}{2}\}) \geq 1 -  \alpha$ is the same. 

Due to symmetry, we have, for any $j = 1,\ldots,n$,
\begin{equation}
\label{eq:A1}
\mathbb{E}{\bf 1}_{\left\{{\sum_{i \in {\cal B}_j} {\bf 1}_{\{R_{ji} \ge R_{ij} \}}} \ge  \frac{\alpha n}{2} \right\}}= 
\mathbb{E}{\bf 1}_{\big\{{\sum_{i = 1}^n {\bf 1}_{\{R_{new,i} \ge R_{i,new} \}}} \ge \frac{\alpha n}{2}  \big\}},
\end{equation}
where ${\cal B}_j = \{i: (\X_i, y_i) \in {\cal D}_{obs} \cup \{(\X_{new}, y_{n + 1})\}, i \not = j \}$. 
Thus, by the definition of $Q_{n}(y_{new})$ and equation (\ref{eq:A1}), we have 
\begin{align}
   & \mathbb{P}\left(y_{new}\in \left\{y:  Q_{n}(y) \ge \frac{\alpha}{2}\right\}\right)  = \mathbb{E}{\bf 1}_{\big\{{\sum_{i = 1}^n {\bf 1}_{\{R_{new,i} \ge R_{i,new} \}}} \ge \frac{\alpha n}{2}  \big\}} \nonumber \\ 
   & \qquad  = \frac1n
   \sum_{j = 1}^n  \mathbb{E}{\bf 1}_{\left\{{\sum_{i \in {\cal B}_j} {\bf 1}_{\{R_{ji} \ge R_{ij} \}}} \ge  \frac{\alpha n}{2} \right\}} 
    =
    \frac1n \mathbb{E}\{N_n\}
    \label{eq:A2}
\end{align}
where $N_n =
\sum_{j = 1}^n  {\bf 1}_{\left\{{\sum_{i \in {\cal B}_j} {\bf 1}_{\{R_{ji} \ge R_{ij} \}}} \ge  \frac{\alpha n}{2} \right\}}$
is the
size (number of elements) of the set 
${\cal J} = \{j'| $ $\sum_{i \in {\cal B}_{j'}} {\bf 1}_{\{R_{j'i} \ge R_{ij'} \}} \ge \frac{\alpha n}{2},j' = 1 \ldots, n \}$.

To give an lower bound on $N_n$, consider for any $j\notin \cal J$, we have
\begin{align*}
    \frac{\alpha n}{2} &> \sum_{i \in {\cal B}_j} {\bf 1}_{\{R_{ji} \ge R_{ij} \}} = \sum_{i \in  {\cal B}_j \cap{\cal J}} {\bf 1}_{\{R_{ji} \ge R_{ij} \}} + \sum_{i \in {\cal B}_j\backslash \cal J} {\bf 1}_{\{R_{ji} \ge R_{ij} \}} \\ & \geq \sum_{i \in {\cal B}_j\backslash \cal J} {\bf 1}_{\{R_{ji} \ge R_{ij} \}}
\end{align*}
Summing over all $j\notin \cal J$ (which has $n - N_n$ members) and by symmetry, we have
\begin{align*}
     \frac{\alpha n}{2}(n-N_n) & > \sum_{j \not \in \cal J}\sum_{i \in {\cal B}_j\backslash \cal J} {\bf 1}_{\{R_{ji} \ge R_{ij} \}} = \sum_{j \not \in {\cal J}, i \not \in { \cal J}, i \not = j}{\bf 1}_{\{R_{ji} \ge R_{ij} \}} \\ & = \frac{(n-N_n)(n-N_n-1)}{2}
\end{align*}
Solving the above inequality, we get a lower bound on $N_n\geq (1-\alpha)n$. By (\ref{eq:A2}), it follows that $\mathbb{P}(y_{new}\in \{y: Q_n(y) \ge \frac{\alpha}{2}\}) \geq 1 -  \alpha$.

\subsection*{S.2 \, Proof of Proposition 2}

{\bf Lemma A1:} Suppose that  $\X_i$ are IID samples from a distribution ${\cal F}_x$ and $ \mathbb{E}(h_{ii})$ and $ \mathbb{E}(g_{ii})$ exist, for $i = 1, \ldots, n$. Then, for any fixed $0<\kappa <1$, $h_{ii} = o_p(n^{-\kappa})$ and $g_{ii} = o_p(n^{-\kappa})$, $i = 1, \ldots, n$.

{\bf Proof of Lemma A1.} Due to IID symmetric, we have $ \mathbb{E}(h_{11}) = \ldots =  \mathbb{E}(h_{nn}) =  \mathbb{E} (\frac{1}{n}\sum_{i = 1}^n h_{ii})$ $= \frac{p}{n}$. By Markov inequality, $P\big(n^{\kappa} h_{ii}>\delta \big)\leq \frac{p/n}{\delta/n^\kappa}\to 0$, for any $\delta>0$. Thus, $h_{ii} = o_p(n^{-\kappa})$. The proof of $g_{ii} = o_p(n^{-\kappa})$ is the same. 

{\bf Lemma A2:}
Suppose $\X_{new} = \bar \X$ and $\Sigma_{\rm w|z} ={\rm var}\left(\mathbf{w}_i|\mathbf{z}_i\right)> 0$. Define
$U_i^{(n)} = (1-h_{i,new})u_i$ and $V_i^{(n)}=$ $ (1-g_{i,new})v_i$, $i = 1, \ldots, n$.  Then, as $n \to \infty$,
{\small
\begin{equation}
\label{eq:quantile1}
\frac1n \mathbb{E}
\left\{\big|\sum_{i=1}^n{\bf 1}_{\{U_i^{(n)}\leq x\}} - \sum_{i=1}^n\frac12 \big({\bf 1}_{\{\epsilon_i\leq x+\tau_n\}}+{\bf 1}_{\{\epsilon_i\leq x-\tau_n\}}\big)\big|  \bigg| {\bf X} \right\} \to 0
\end{equation}
}
\vspace{-4mm}
\noindent
and
{\small
\begin{equation}
\label{eq:quantile2}
\frac1n\mathbb{E}\left\{\big|\sum_{i=1}^n{\bf 1}_{\{V_i^{(n)}\leq x\}} - \sum_{i=1}^n\frac12 \big({\bf 1}_{\{\epsilon_i+w_i^T\beta_2\leq x+\tau_n\}}+{\bf 1}_{\{\epsilon_i+w_i^T\beta_2\leq x-\tau_n\}}\big)\big| \bigg| {\bf Z} \right\} \to 0,
\end{equation}
}
\hspace{-5mm}

\noindent
for any given $\tau_n = n^{-\frac{\kappa}{2}}$, $0<\kappa<1$.

{\bf Proof of Lemma A2.}
Since $\X_{new} = \bar \X$, we have $h_{i,new} = \frac1n\sum_{j=1}^n h_{ij}=\frac1n$~and  $g_{i,new}=\frac1n\sum_{j=1}^n g_{ij} = \frac1n$. Thus, $U_i^{(n)} = \frac{1-1/n}{1-h_{ii}}\big(\epsilon_i-\sum_{j=1}^n h_{ij}\epsilon_j\big)$ and $V_i^{(n)} = \frac{1-1/n}{1-g_{ii}}\big(e_i-\sum_{j=1}^n g_{ij}e_j\big)$, where $e_i = \epsilon_i+w_i^T\beta_2$.  
Furthermore, since $\sum_{j=1}^n h_{ij}^2 = h_{ii}$ and $\sum_{j=1}^n g_{ij}^2 = g_{ii}$, it follows that, given $\mathbf{X}$, the conditional distribution of  $\sum_{j=1}^n h_{ij}\epsilon_j \sim  N(0,\sigma^2h_{ii})$ and, given $\mathbf{Z}$, the conditional distribution of
$\sum_{j=1}^n g_{ij}e_j\sim N(0,g_{ii}(\sigma^2 + \beta_2^T \Sigma_{\rm w|z}\beta_2))$. 
We have, as $n \to \infty$,

\vspace{-4mm}
{\small
\begin{align*}
 &\frac1n\mathbb{E}\left\{\big|\sum_{i=1}^n{\bf 1}_{\{U_i^{(n)}\leq x\}} -\sum_{i=1}^n\frac12\big\{{\bf1}_{\{\epsilon_i\leq x+\tau_n\}}+ {\bf 1}_{\{\epsilon_i\leq x-\tau_n\}}\big\}\big|\bigg|\mathbf{X}\right\}\nonumber\\
    &\leq \frac1n\sum_{i=1}^n   \mathbb{E}\left\{\big|{\bf 1}_{\{U_i^{(n)}\leq x\}} - \frac12({\bf 1}_{\{\epsilon_i\leq x+\tau_n\}}+{\bf 1}_{\{\epsilon_i\leq x-\tau_n\}})\big|\bigg|\mathbf{X} \right\}\nonumber\\
    &= \frac1n\sum_{i=1}^n\bigg[ \mathbb{E}\left({\bf 1}_{\{U_i^{(n)}\leq x,\epsilon_i> x+\tau_n\}}\bigg|\mathbf{X} \right) + \mathbb{E}\left({\bf 1}_{\{U_i^{(n)}> x,\epsilon_i\leq x-\tau_n\}}\bigg|\mathbf{X} \right)\\
   &\quad\quad + \frac12\mathbb{E}\left({\bf 1}_{\{x-\tau_n<\epsilon_i\leq x+\tau_n\}}({\bf 1}_{\{U_I^{(n)}\leq x\}}+{\bf 1}_{\{U_i^{(n)}>x\}})\bigg|\mathbf{X} \right)\bigg] \nonumber\\
    &= \frac1n\sum_{i=1}^n\bigg[ \mathbb{P}\bigg(\frac{1-\frac1n}{1-h_{ii}}\big(\epsilon_i-\sum_{j=1}^n h_{ij}\epsilon_j\big)\leq x,\epsilon_i> x+\tau_n\big|\mathbf{X} \bigg) \\ &\quad\quad + \mathbb{P}\bigg(\frac{1-\frac1n}{1-h_{ii}}\big(\epsilon_i-\sum_{j=1}^n h_{ij}\epsilon_j\big)> x,\epsilon_i\leq x-\tau_n\big|\mathbf{X} \bigg)
    + \frac12\mathbb{P}\bigg(x-\tau_n<\epsilon_i\leq x+\tau_n\big|\mathbf{X} \bigg)\bigg] \nonumber
    \\
    &\leq \frac1n\sum_{i=1}^n\bigg[ \mathbb{P}\bigg(\sum_{j=1}^n h_{ij}\epsilon_j> \frac{h_{ii}-\frac1n}{1-\frac1n}x+\tau_n\big|\mathbf{X} \bigg)  + \mathbb{P}\bigg(\sum_{j=1}^n h_{ij}\epsilon_j< \frac{h_{ii}-\frac1n}{1-\frac1n}x-\tau_n\big|\mathbf{X} \bigg) \\
    &\quad\quad 
    + \frac12\left\{ \Phi\big(\frac{x+\tau_n}{\sigma})- \Phi(\frac{x-\tau_n}{\sigma}\big)\right\}\bigg] \\
    &\leq \frac1n \sum_{i=1}^n \mathbb{P}\bigg( \big|\sum_{j=1}^n h_{ij}\epsilon_j\big|> \frac{\tau_n}{2}\bigg|\mathbf{X}  \bigg)+\frac{\tau_n}{2 \pi \sigma^2}\leq \frac1n \sum_{i=1}^n \frac{4 h_{ii}\sigma^2 }{\tau_n^2}+\frac{\tau_n}{2 \pi \sigma^2} = \frac{4p \sigma^2}{(n \tau_n)^2}+  \frac{\tau_n}{2 \pi \sigma^2}.
\end{align*}
}

\vspace{-3.4mm}
\noindent
The second last inequality holds since $\tau_n = n^{-\frac{\kappa}2}$, $0 < \kappa < 1$, and by Lemma A1, $\frac{h_{ii} |x|}{\tau_n^2} < \frac 12$ for $n$ large enough. 
Thus, (\ref{eq:quantile1}) follows. Similarly, we can prove (\ref{eq:quantile2}).

{\bf Proof of Proposition 2.}
Denote by $h_n(x) = \frac1n\sum_{i=1}^n{\bf 1}_{\{U_i^{(n)}\leq x\}}$ and $g_n(x) = \frac1n\sum_{i=1}^n{\bf 1}_{\{\epsilon_i\leq x\}} $. 
Write $a_\alpha = q_\alpha(\{U_i^{(n)}\}_{i=1}^n )$ and $b_\alpha = q_\alpha(\{\epsilon_i\}_{i=1}^n )$, for $\alpha\in (0,1)$, and they are solutions to equation
$h_n(x) = \alpha$ and $g_n(x) = \alpha$, respectively. Also, define $c_\alpha$ to be the solution of equation $\frac12 \big\{g_n(x-\tau_n)+g_n(x+\tau_n)\big\} = \alpha$. 

Since $g_n(x)$ increases in $x$, 
$g_n(x-\tau_n)\leq \frac12 \big\{g_n(x-\tau_n)+g_n(x+\tau_n)\big\}\leq g_n(x+\tau_n)$. It follows that $b_\alpha - \tau_n \leq c_\alpha \leq b_\alpha + \tau_n$. That is,  $|c_\alpha - b_\alpha| = o_p(1)$. Since $\epsilon$ is independent of $\mathbf{X}$, this statement also holds conditional on $\mathbf{X}$.
On the other hand, by Lemma A2, conditional on $\mathbf{X}$, we have $|h_n(x) - \frac12 \big\{g_n(x-\tau_n)+g_n(x+\tau_n)\big\}| = o_p(1)$, for $x \in (- \infty, \infty)$. Since both $h_n(x)$ and $g_n(x)$ are increasing functions in $x$, 
it follows that, conditional on $\mathbf{X}$, $|a_{\alpha} - c_{\alpha}| = o_p(1)$. 
Putting things together, we have $|a_{\alpha} - b_{\alpha}| = o_p(1)$; That is, conditional on~$\mathbf{X}$,
\begin{equation}
    \label{eq:AA}
|q_\alpha(\{U_i^{(n)}\}_{i=1}^n ) -  q_\alpha(\{\epsilon_i\}_{i=1}^n )| = o_p(1).
\end{equation}
It follows that, conditional on $\mathbf{X}$ (thus conditional on $\mathbf{Z}$
 and unconditionally), 
{\small \begin{equation}\label{eq:LengthU}
q_{1-\frac{\alpha}2}\left(\{U_i^{(n)}\}_{i=1}^n \right) - q_{\frac{\alpha}2}\left(\{U_i^{(n)}\}_{i=1}^n \right)
    = \sigma \left\{\Phi^{-1}\left( 1- \frac{\alpha}{2}\right)
     - \Phi^{-1}\left( \frac{\alpha}{2}\right)\right\}
    + o_p(1). 
\end{equation}
}

\vspace{-4mm}
Since $\X_i = (\Z_i^T, \W_i^T)^T$ are independent copies and $\epsilon_i$ is independent of  $\X_i$, we have,     conditional on  $\mathbf{Z}$, $\epsilon_i+w_i^T\beta_2 \overset{iid}{\sim} N(0, \sigma_{\Sigma}^2)$, where $\sigma_{\Sigma}^2 = \sigma^2+\beta_2^T\Sigma_{\rm w|z}\beta_2$.
Similar to the proof of (\ref{eq:AA}), we can show that, conditional on $\mathbf{Z}$,  
$$
\big| q_{\alpha}\big(\{V_i^{(n)}\}_{i=1}^n\big) - q_{\alpha}\big(\{\epsilon_i+w_i^T\beta_2\}_{i=1}^n\big)\big| = o_p(1). 
$$
It follows immediately that, conditional on  $\mathbf{Z}$ (thus also unconditionally),
{\small
\begin{equation}\label{eq:LengthV}
q_{1-\frac{\alpha}2}\left(\{V_i^{(n)}\}_{i=1}^n \right) - q_{\frac{\alpha}2}\left(\{V_i^{(n)}\}_{i=1}^n \right)
    = \sigma_{\Sigma} \left\{\Phi^{-1}\left( 1- \frac{\alpha}{2}\right)
     - \Phi^{-1}\left( \frac{\alpha}{2}\right)\right\}
    + o_p(1). 
\end{equation}
}

\vspace{-4mm}
Finally, since $\sigma_{\Sigma}^2 > \sigma^2$, Proposition 2 holds by (\ref{eq:LengthU})~and~(\ref{eq:LengthV}).

\end{document}